\documentclass[11pt]{article}

\usepackage{amssymb,amsthm,amsmath,hyperref,txfonts}
\usepackage{color}

\usepackage{graphicx,float}

\numberwithin{equation}{section}
\newtheorem{thm}{Theorem}[section]
\newtheorem{lem}{Lemma}[section]

\newcommand{\beq}{\begin{eqnarray}}
\newcommand{\eeq}{\end{eqnarray}}
\newcommand{\beqno}{\begin{eqnarray*}}
\newcommand{\eeqno}{\end{eqnarray*}}
\newcommand{\be}{\begin{equation}}
\newcommand{\ee}{\end{equation}}
\theoremstyle{definition}
\newtheorem{defn}{Definition}[section]

\newcommand{\D}{\displaystyle}

\allowdisplaybreaks

\begin{document}

\title{\bf Energy conservation for the weak solutions to the equations of compressible magnetohydrodynamic flows in three dimensions}
\author{Tingsheng Wang\thanks{E-mail: wangtingsheng365@163.com.},\ \
 Xinhua Zhao\thanks{E-mail:
xinhuazhao211@163.com. Corresponding author.},\ \ Yingshan Chen\thanks{Email: mayshchen@scut.edu.cn.},\ \ Mei Zhang\thanks{Email: scmzh@scut.edu.cn.}
\\
\textit{\small  School of Mathematics, South China
University of Technology, Guangzhou 510641, China} \\
}
\date{\today}
\maketitle
\begin{abstract}
In this paper, we prove the energy conservation for the weak solutions to the three-dimensional equations of compressible magnetohydrodynamic flows (MHD) under certain conditions only about density and velocity. This work is inspired by the seminal work by Yu \cite{Yu2017} on the energy conservation of compressible Navier-Stokes equations. Our result indicates that even the magnetic field is taken into account, we only need
some regularity conditions of the density and velocity as in \cite{Yu2017} to ensure the energy conservation.

\end{abstract}
\vspace{4mm}



\section{Introduction And Main Results}

%~~~~~~Magnetohydrodynamics (MHD) is the study of the interaction between magnetic fields and moving, conducting fluids.
% Example of such magnetofluids include plasmas, liquid metals, salt water, and electrolytes.
% The field of MHD was initiated by Hannes Alfv\'{e}n in \cite{Alfven1} to explain energy transport in the Sun,
%and the field of plasma astrophysics quickly evolved \cite{Alfven}.
%for which he received the Nobel Prize in Physics in 1970. Moreover, geophysical questions such as the maintenance of the Earth's magnetic field and
%the interaction of the solar wind with the magnetosphere continue to be addressed \cite{Cowling,Davidson}.
Magnetohydrodynamics (MHD) concerns the motion of conducting fluids in
an electromagnetic field and has a very broad range of applications. The dynamic
motion of the fluid and the magnetic field interact strongly on each other. In this paper, the fluid we consider is isentropic and compressible, namely, it is governed by the isentropic compressible Navier-Stokes equations. The equations of the magnetic field are called the induction equation.%The set of equations that describe MHD are a combination of the %Navier-Stokes equations of fluid dynamics
%and Maxwell's equations of electro-magnetism. In this paper we are interested in the three-dimensional
%compressible magnetohydrodynamic flows in the isentropic case,
Hence the compressible MHD system for isentropic flows can be written as below \cite{Cabannes,Kulikovskiy,Laudau}.
\begin{equation}\label{2dbu-E1.1}
\left\{
\begin{array}{l}
\rho_{t} + \mathrm{div}(\rho \mathbf{u}) = 0,\\
(\rho \mathbf{u})_t + \mathrm{div}(\rho \mathbf{u}\otimes \mathbf{u}) + \nabla P = (\nabla \times \mathbf{H}) \times \mathbf{H} + \mu \Delta\mathbf{u} + (\lambda + \mu)\nabla(\mathrm{div}\mathbf{u}),\\
\mathbf{H}_t - \nabla \times (\mathbf{u} \times \mathbf{H}) = -\nabla \times(\nu \nabla \times \mathbf{H}), \ \mathrm{div}\mathbf{H} = 0,\\
\end{array}
\right.
\end{equation}
where $\rho=\rho(x,t), \mathbf{u}(x,t) = (u_1, u_2, u_3)(x,t), \mathbf{H} = (H_1,H_2,H_3)(x,t)$ denote the density of the fluid, the velocity field
and the magnetic field, respectively; $P(\rho) = a\rho^{\gamma}$ is the pressure with constants $a>0$, and $\gamma>1$;
the constants $\mu$ and $\lambda$ are the shear and bulk viscosity coefficients satisfying
the physical restriction $3\lambda + 2\mu \ge0$ and $\mu>0$; and the constant $\nu>0$ is the magnetic diffusivity.
The positive constant $a$ does not play essential role in the following analysis. Thus for simplicity we take $a=1$.
%and all these kinetic coefficients and the magnetic diffusivity are independent of the magnitude and direction of the magnetic field.
%The symbol $\otimes$ denotes the Kronecker product. Usually, we refer to the first equation in \eqref{2dbu-E1.1} as the continuity equation,
%the second equation as the momentum balance equation, and the last equation in \eqref{2dbu-E1.1} is called the induction equation.\\

For the sake of simplicity we
will consider the case of a bounded domain with periodic boundary conditions in $\mathbb{R}^3$, namely $\Omega=\mathbb{T}^3$,
and the following initial conditions:
\begin{equation}\label{I.C.}
(\rho, \rho\mathbf{u}, \mathbf{H})(x,0) = (\rho_0, \mathbf{m}_0, \mathbf{H}_0)(x), ~~x\in\Omega,
\end{equation}
where we define $\mathbf{m}_0=\mathbf{0},$ if $\rho_0=0.$

%In this paper, we study the multi-dimensional isentropic problem (1.1) and (1.2) with $\gamma > \frac32$, where all the viscosity coefficients $\mu, \lambda,\nu$ are constant.

%In this paper, we study the multi-dimensional isentropic problem (1.1) and (1.2) with $\gamma > \frac32$, where all the viscosity coefficients $\mu, \lambda,\nu$ are constant.

The global existence of weak solutions to (\ref{2dbu-E1.1}) in a bounded domain of $\mathbb{R}^3$ was obtained by Hu and Wang \cite{Wang2010} for $\gamma>\frac{3}{2}$. Moreover, the global weak solutions exist in the renormalized sense with arbitrarily large initial data as well, satisfying the energy inequality
\beq \label{energy inequ}
&&\int_\Omega\bigg(\frac12 \rho \mathbf{u}^2 +\frac{\rho^\gamma }{\gamma - 1} + \frac12|\mathbf{H}|^2\bigg)\,\mathrm{d}x+ \int_0^t \int_\Omega\bigg(\mu|\nabla\mathbf{u}|^2 + (\lambda+\mu)(\mathrm{div}\mathbf{u})^2+\nu|\nabla\times\mathbf{H}|^2\bigg)\,\mathrm{d}x\,\mathrm{d}t\notag\\
&&\le\int_\Omega\bigg(\frac12 \rho_0 \mathbf{u}_0^2 +\frac{\rho_0^\gamma }{\gamma - 1} + \frac12 |\mathbf{H}_0|^2\bigg)\,\mathrm{d}x
\eeq for $t\in(0,\infty)$. In fact, when the solutions are smooth enough such as strong solutions or classical solutions,
the energy inequality (\ref{energy inequ-1}) can be written as an equality, namely,
\beq\label{energy inequ-1}
&&\int_\Omega\bigg(\frac12 \rho \mathbf{u}^2 +\frac{\rho^\gamma }{\gamma - 1} + \frac12|\mathbf{H}|^2\bigg)\,\mathrm{d}x+ \int_0^t \int_\Omega\bigg(\mu|\nabla\mathbf{u}|^2 + (\lambda+\mu)(\mathrm{div}\mathbf{u})^2+\nu|\nabla\times\mathbf{H}|^2\bigg)\,\mathrm{d}x\,\mathrm{d}t\notag\\
&&=\int_\Omega\bigg(\frac12 \rho_0 \mathbf{u}_0^2 +\frac{\rho_0^\gamma }{\gamma - 1} + \frac12 |\mathbf{H}_0|^2\bigg)\,\mathrm{d}x
  \eeq for $t\in(0,\infty)$. For example, see
%There are a lot of studies on isentropic compressible MHD equations by many physicists and mathematicians,
%due to its physical importance, complexity, rich phenomena, and mathematical challenges, such as %\cite{CW1,CW2,Fan,hu2008,Hu2008,Wang2010,Kawashima,Li2013,Volpert} and the references therein.
%For the incompressible isentropic MHD equations (where $\rho$ is a constant), Duraut and Lions \cite{DuvautandLions} constructed a class of a global weak solutions with finite energy. For the three-dimensional compressible MHD equations, Umeda-Kawas-Shizuta \cite{Umeda} obtained the global existence of smooth solutions to the linearized MHD equations.
  \cite{CW1, CW2, Kawashima, Volpert, Li2013} for global smooth solutions in one dimension with arbitrarily large initial data and in multi-dimensions with small perturbations of a given constant
state, and \cite{Fan, Volpert} for local strong solutions with arbitrarily large initial data.

The question is how much regularity of the weak solutions is needed to ensure the energy equality (\ref{energy inequ-1})?
%The local strong solution to the compressible MHD equations with general initial data was proved by Volpert and Khudiaev \cite{Volpert}. Recently, Hu and Wang \cite{HuWang1,HuWang2,Wang2010} established the existence of global weak solutions to the compressible MHD equations with general initial data by the Faedo-Galerkin method and the vanishing viscosity method. However due to the strong interaction between the dynamic motion of fluid and magnetic field, both hydrodynamic and electrodynamic are coupled together, which makes the MHD problem difficult. For example, even for the one-dimensional case, the global existence of a classical solution to the full perfect MHD equations with large data remains unsolved when all the viscosity, heat conductivity, and diffusivity coefficients was constant, although the corresponding problem for the Navier-Stokes equations was solved in \cite{Kazhikhov}.
In the context of incompressible Euler equations, this
question is linked to a famous conjecture of Onsager \cite{Onsager}. %More precisely, Onsager conjectured that, in the context of homogeneous
%incompressible Euler equations, kinetic energy is globally conserved for H\"{o}lder continuous solutions with
%the exponent greater than 1/3, while energy dissipation phenomenon occurs for H\"{o}lder continuous solutions
%with the exponent less than 1/3.
It has been made great progress recently \cite{Bardos2018, Buckmaster2015, Buckmaster2016, Cheskidov, Constantin, Eyink}. %Latter, Isett \cite{Isett} solved the second part of the
%Onsager's conjecture.
In the context of incompressible Navier-Stokes equations, Serrin \cite{Serrin} proved the energy conservation under the condition $\mathbf{u}\in L^p(0,T;L^q(\Omega))$,
$\frac2p+\frac Nq\leq1$, where N is the dimension. Later, Shinbrot \cite{Shinbrot} removed
the dimensional dependence, i.e., $\frac2p+\frac 2q\leq1$, where $q\geq4$. When the magnetic field is ignored, i.e. $\mathbf{H}=0$, system \eqref{2dbu-E1.1} becomes the compressible Navier-Stokes equations.  Yu \cite{Yu2017} proved the energy conservation (\ref{energy inequ-1}) ($\mathbf{H}=0$) of the Lions-Feireisl weak solutions (see \cite{Feireisl,Feir2,Lions}) for $\Omega=\mathbb{T}^3$ provided that
 \beq\begin{cases}
0\le\rho\le \tilde{\rho}<\infty,\, \mathrm{and}\, \, \nabla\sqrt{\rho}\in L^\infty\big(0,T;L^2(\Omega)\big),\\[2mm]u\in L^p(0,T;L^q(\Omega)),\, \, \mathrm{for\, any}\,\, \frac{1}{p}+\frac{1}{q}\leq\frac{5}{12},\, \mathrm{and}~q\geq6,\\[2mm]
u_0\in L^k(\Omega),\, \frac{1}{k}+\frac{1}{q}\le\frac{1}{2},
 \end{cases}\eeq where $\tilde{\rho}$ is a positive constant. In \cite{Yu2017}, the case of density-dependence viscosity is also considered.
Recently, Chen, Liang, Wang, Xu \cite{Wang2018} nicely extended Yu's results to the Dirichlet problem.
%The method of mollification and estimation of commutator errors was employed to prove that, if a weak solution $\mathbf{u}$ of the incompressible Euler system belongs to $L^3([0,T],B^{\alpha,\infty}_3(\mathbb{T}^3))$. Roughly speaking, Onsager's conjecture (for Euler equations) addresses the equation how much regularities needed for a weak solution to conserve energy.  In the context of classical incompressible Navier-Stokes equations, since global regularity in three dimensions has been a long standing open question, it is natural to ask how much regularities needed for a weak solution to satisfy the energy

The purpose of this paper is to provide a sufficient condition for the energy conservation of the weak solution of
\eqref{2dbu-E1.1}-\eqref{I.C.}, which is motivated by Yu's work \cite{Yu2017} (see also \cite{Wang2018}). %If the solutions of \eqref{2dbu-E1.1}-\eqref{I.C.} is sufficiently smooth, the energy conservation holds. However, the existence of smooth solutions is
%a longstanding open question. So this question could be transformed to ask how badly of  a weak solution can keep the
%energy conservation. Mathematically, what is minimal regularity such that a weak solution satisfies the following energy equality
%\begin{equation}\begin{split}
%\int_\Omega\bigg(\frac12 \rho \mathbf{u}^2 &+\frac{\rho^\gamma }{\gamma - 1} + \frac12|\mathbf{H}|^2\bigg)\,\mathrm{d}x+ \int_0^T %\int_\Omega\bigg(\mu|\nabla\mathbf{u}|^2 + %(\lambda+\mu)(\mathrm{div}\mathbf{u})^2+\nu|\nabla\times\mathbf{H}|^2\bigg)\,\mathrm{d}x\,\mathrm{d}t\notag\\
%&=\int_\Omega\bigg(\frac12 \rho_0 \mathbf{u}_0^2 +\frac{\rho_0^\gamma }{\gamma - 1} + \frac12 |\mathbf{H}_0|^2\bigg)\,\mathrm{d}x.
%\end{split} \end{equation}

%It is well known that if a solution is smooth enough, then it conserves the energy.
%Thus, our question is connected to the regularity of weak solutions.

%Multiplying the first equation in (1.1) by $b'(\rho)$, we obtain the renormalized continuity equation:
%\begin{equation}\label{rce}
%b(\rho)_t+\mathrm{div}(b(\rho)\mathbf{u})+(b'(\rho)-b(\rho))\mathrm{div}\mathbf{u}=0,
%\end{equation}
%for some suitable function $b \in C^1(\mathbb{R}^{+})$.

%Following the strategy in \cite{Feireisl,Lions,Vasseur}, we introduce the concept of weak solution $(\rho, \mathbf{u}, \mathbf{H})$ to %\eqref{2dbu-E1.1}-\eqref{I.C.}.
\begin{defn}\label{defn} (weak solution)
$(\rho,\mathbf{u},\mathbf{H})$ is called a weak solution to \eqref{2dbu-E1.1}-\eqref{I.C.} over $\Omega\times(0,T)$,
if $(\rho,\mathbf{u},\mathbf{H})$ satisfies that\\
$\bullet$ \eqref{2dbu-E1.1} holds in $\mathcal{D}^{\prime}\big(\Omega \times(0, T)\big)$ satisfying\\
\begin{equation}\label{regu}
 \begin{cases}
 \rho\in L^{\infty}\big(0,T;L^{\gamma}(\Omega)), \quad \rho\geq0,\\
 \sqrt{\rho}\mathbf{u}\in L^\infty(0,T;L^2(\Omega)), \quad \nabla \mathbf{u}\in L^2(0,T;L^2(\Omega)),\\
 \mathbf{H}\in L^\infty(0,T;L^2(\Omega)), \quad \nabla \mathbf{H}\in L^2(0,T;L^2(\Omega)),
 \end{cases}
\end{equation}
and \begin{equation}\label{rho u}
 \rho \mathbf{u}\in C([0,T];L_{weak}^2(\Omega)),
\end{equation}
and
\begin{equation}\label{H}
  \mathbf{H}\in C([0,T];L_{weak}^2(\Omega));
\end{equation}
$\bullet$ the energy inequality (\ref{energy inequ}) holds;

\noindent$\bullet$ \eqref{I.C.} holds in $\mathcal{D}^{\prime}(\Omega)$.
\end{defn}

\medskip

Our main result reads as follows.
\begin{thm}\label{thm}
Assume that
\begin{equation}\label{initial}
 \begin{cases}
 \rho(x, 0)=\rho_{0}(x)  \in  L^1(\Omega)\cap L^{\gamma}(\Omega), \quad \rho_{0}(x) \geq 0~a.e.~in~\Omega,\\
 \rho(x, 0) \mathbf{u}(x, 0)=\mathbf{m}_{0}(x) \in L^{1}(\Omega),  \quad \mathbf{m}_{0}=\mathbf{0} \text { if } \rho_{0}=0, \quad \frac{\left|\mathrm{m}_{0}\right|^2}{\rho_{0}} \in L^{1}(\Omega),\\
  \mathbf{H}(x, 0)=\mathbf{H}_{0}(x) \in L^{2}(\Omega), \quad \operatorname{div} \mathbf{H}_{0}=0 \text { in } \mathcal{D}^{\prime}(\Omega).
 \end{cases}
\end{equation}
In addition, we assume $\mathbf{u}_0\in L^\kappa,$ where $\kappa>2$. Let $(\rho,\mathbf{u},\mathbf{H})$ be a weak solution to \eqref{2dbu-E1.1}-\eqref{I.C.} in the sense of Definition \ref{defn}.
Moreover, if
\begin{equation}\begin{split}\label{rho}
0\leq\rho\leq\tilde{\rho}<\infty,~~~\nabla\sqrt{\rho}\in L^{\infty}(0,T;L^2(\Omega)).
\end{split}\end{equation}
and
\begin{equation}\label{u}
\mathbf{u}\in L^p(0,T;L^q(\Omega)) ~~~~for ~any~\frac{2}{p}+\frac{3}{q}\leq1,~ with ~q\geq6,
\end{equation}
then the weak solution $(\rho,\mathbf{u},\mathbf{H})$ satisfies the energy equality (\ref{energy inequ-1}) for $t\in[0,T]$.
%\begin{equation}\begin{split}
%\int_\Omega\bigg(\frac12 \rho \mathbf{u}^2 &+\frac{\rho^\gamma }{\gamma - 1} + \frac12|\mathbf{H}|^2\bigg)\,\mathrm{d}x+ \int_0^T %\int_\Omega\bigg(\mu|\nabla\mathbf{u}|^2 + (\lambda+\mu)(\mathrm{div}\mathbf{u})^2+\nu|\nabla\times\mathbf{H}|^2\bigg)\,\mathrm{d}x\,\mathrm{d}t\\
%&=\int_\Omega\bigg(\frac12 \rho_0 \mathbf{u}_0^2 +\frac{\rho_0^\gamma }{\gamma - 1} + \frac12 |\mathbf{H}_0|^2\bigg)\,\mathrm{d}x.
%\end{split}\end{equation}
\end{thm}

\section{Preliminaries}

%~~~~To state the main results, let us introduce some notations and inequalities.\\
Define
\begin{align*}
\overline{f(x,t)} &= \eta_{\epsilon}*f(x,t)
= \int_0^t \int_\Omega\eta_{\epsilon}(x-y,t-s)f(y,s)\,\mathrm{d}y\,\mathrm{d}s\\
&= \int_0^t \int_\Omega\frac{1}{\epsilon^{4}}\eta(x-y,t-s)f(y,s)\,\mathrm{d}y\,\mathrm{d}s
\end{align*}
where $\eta_{\epsilon}(x,t) = \frac{1}{\epsilon^4}\eta(\frac{x}{\epsilon},\frac{t}{\epsilon})$, and $\eta(t,x)\geq0$ is a smooth even function compactly supported in the space-time ball of radius 1, and with an integral equal to 1.

%Note that, for any weak solution $(\rho, \mathbf{u}, \mathbf{H})$, with condition \eqref{rho}, Theorem \ref{thm} satisfies:
%\begin{align}\label{regularity of u,H}
%\|\mathbf{H}\|_{L^{\infty}(0,T;L^2(\Omega))}\leq C<\infty,\ \|\nabla\mathbf{H}\|_{L^2(0,T;L^2{\Omega})}\leq C<\infty.
%\end{align}

%\begin{lem}\label{domain}
%For any weak solution $(\rho,\mathbf{u})$ in the sense of Definition 1.1, with additional conditions, $\rho_t$ is bounded in $L^p(0,T;L^\frac{2q}{q+2}(\Omega)) + L^2(0,T;L^2(\Omega))$. In particular, $\rho_t$ is bounded in $L^2(0,T;L^{\frac{2q}{q+2}}(\Omega))$ if $p\geq2$.
%\end{lem}
The following lemma will be useful in the proof of Theorem \ref{thm}.
%\begin{lem}[Lemma 2.2, \cite{Yu2017}]\label{inequa}
%If $\mathbf{u}\in L^{\infty}(0,T;L^2(\Omega))\cap L^p(0,T;L^q(\Omega))$, then there exists some $\alpha_1\in(0,1)$ such that $\mathbf{u}\in L^{r}(0,T;L^s(\Omega))$, for any
%\begin{align*}
%\frac{1}{r}=\frac{1-\alpha_1}{p},
%\end{align*}
%and
%\begin{align*}
%\frac1s=\frac{\alpha_1}{2}+\frac{1-\alpha_1}{q}.
%\end{align*}
%\end{lem}

\begin{lem}[\cite{Lions}]\label{limit}
Let $\partial$ be a partial derivative in space or time. Let $f, \partial f\in L^p(\Omega\times\mathbb{R}^+), g\in L^q(\Omega\times\mathbb{R}^+)$ with $1\leq p,q\leq\infty$, and $\frac1p+\frac1q \leq 1$. Then, we have
\begin{align*}
\|\overline{\partial(fg)}-\partial(f\overline{g})\|_{L^r(\Omega\times\mathbb{R}^+)}\leq C\|\partial f\|_{L^p(\Omega\times\mathbb{R}^+)}\|g\|_{L^q(\Omega\times\mathbb{R}^+)}
\end{align*}
for some constant $C>0$ independent of $\epsilon$, f and g, and with $\frac1r = \frac1p + \frac1q$. In addition,
\begin{align*}
\overline{\partial(fg)}-\partial(f\overline{g})\to 0 \  in\  L^r(\Omega\times\mathbb{R}^+)
\end{align*}
as $\epsilon \to 0$, if $r < \infty$.
\end{lem}

%\begin{lem}
%Let $f\in W^{1,p}(\mathbb{R}^N)$, $g\in L^q(\mathbb{R}^N)$ with $1\leq p,q \leq \infty$, and $\frac1p +\frac1q \leq 1$. Then, we have
%\begin{align*}
%\|\overline{\mathrm{div}(fg)} - \mathrm{div}(f\overline{g})\|_{L^r(\mathbb{R}^N)}\leq C\|f\|_{W^{1,p}(\mathbb{R}^N)}\|g\|_{L^q(\mathbb{R}^N)}
%\end{align*}
%for some $C\geq 0$ independent of $\epsilon$ ,f and g, r is determined by $\frac1r = \frac1p + \frac1q$. In addition
%\begin{align*}
%\overline{\mathrm{div}(fg)} - \mathrm{div}(f\overline{g})\to 0\  in \  L^r(\mathbb{R}^N)
%\end{align*}
%as $\epsilon \to 0$ if $r<\infty$.
%\end{lem}
%This lemma includes the following statement:
%\begin{lem}
%Let $f_t \in L^p(0,T),\ g\in L^q(0,T)$ with $1\leq p,q \leq \infty$, and $\frac1p +\frac1q \leq 1$. Then, we have
%\begin{align*}
%\|\overline{(fg)_t} - (f\overline{g})_t\|_{L^r(0,T)}\leq C\|f_t\|_{L^p(0,T)}\|g\|_{L^q(0,T)}
%\end{align*}
%for some $C\geq0$ independent of $\epsilon$, f anf g, r is determined by $\frac1r = \frac1p + \frac1q$. In addition
%\begin{align*}
%\overline{(fg)_t} - (f\overline{g})_t\to 0\ in\  L^r(0,T)
%\end{align*}
%as $\epsilon \to 0$ if $r<\infty$.
%\end{lem}

\section{Proof of Main Result}
%~~~~~~The main object of this section is to prove our main result.

~~~~For a given test function $\psi(t) \in \mathcal{D}(0,+\infty)$, denote
$\Phi = \overline{\psi(t)\overline{\mathbf{u}}}$.
Since $\mathcal{D}(0,+\infty)$ is a class of all smooth compactly supported functions in $(0,+\infty)$,
$\Phi$ is well defined on $(0,+\infty)$ for $\epsilon$ small enough. Finally, we will extend the result for $\psi(t)\in \mathcal{D}(-1,+\infty)$.\\
%Note that, $\psi(t)$ is compactly supported in $(0,+\infty)$.

\textbf{Step 1.} Choosing $\Phi$ as the test function.\\

Using $\Phi$ as the test function of $\eqref{2dbu-E1.1}_1$, one obtains
\begin{equation}\label{r1}
\int_0^T \int_\Omega \Phi \big((\rho\mathbf{u})_t + \mathrm{div}(\rho\mathbf{u} \otimes \mathbf{u}) + \nabla{P} - (\nabla \times \mathbf{H}) \times \mathbf{H} - \mu \Delta{\mathbf{u}} - (\lambda + \mu)\nabla(\mathrm{div}\mathbf{u})\big)\,\mathrm{d}x\,\mathrm{d}t= 0,
\end{equation}
which in turn yields
\begin{equation}\label{r2}
\int_0^T \int_\Omega \psi(t)\overline{\mathbf{u}}\cdot\bigg(\overline{ (\rho\mathbf{u})_t + \mathrm{div}(\rho\mathbf{u} \otimes \mathbf{u}) + \nabla{P} - (\nabla \times \mathbf{H}) \times \mathbf{H} - \mu \Delta{\mathbf{u}} - (\lambda + \mu)\nabla(\mathrm{div}\mathbf{u})}\bigg)\,\mathrm{d}x\,\mathrm{d}t= 0,
\end{equation}
where we used the fact $\eta(-t,-x) = \eta(t,x)$.

The first two terms in \eqref{r2} yield that
\begin{equation} \begin{split}\label{r3}
&\int_0^T \int_\Omega \psi(t)\overline{(\rho \mathbf{u})_t}\cdot\overline{\mathbf{u}}\,\mathrm{d}x\,\mathrm{d}t
+\int_0^T \int_\Omega \psi(t)\overline{\mathrm{div}(\rho\mathbf{u}\otimes\mathbf{u})}\cdot\overline{\mathbf{u}}\,\mathrm{d}x \,\mathrm{d}t    \\
=& \int_0^T \int_\Omega \psi(t)\Big(\overline{(\rho \mathbf{u})_t} - (\rho \overline{\mathbf{u}})_t\Big)\cdot\overline{\mathbf{u}}\,\mathrm{d}x\,\mathrm{d}t
+\int_0^T \int_\Omega \psi(t) \Big(\overline{\mathrm{div}(\rho\mathbf{u}\otimes\mathbf{u})} - \mathrm{div}(\rho\mathbf{u}\otimes\overline{\mathbf{u}})\Big)\cdot\overline{\mathbf{u}}\,\mathrm{d}x \,\mathrm{d}t\\
&+ \int_0^T \int_\Omega \psi(t)(\rho \overline{\mathbf{u}})_t\cdot\overline{\mathbf{u}}\,\mathrm{d}x \,\mathrm{d}t
+\int_0^T \int_\Omega \psi(t)\Big((\rho\mathbf{u}\cdot\nabla)\overline{\mathbf{u}} + \mathrm{div}(\rho\mathbf{u})\overline{\mathbf{u}}\Big)\cdot \overline{\mathbf{u}}\,\mathrm{d}x \,\mathrm{d}t \\
=&A + B+\int_0^T \int_\Omega \psi(t)\bigg(\frac{1}{2} \rho\mid\overline{\mathbf{u}}\mid^2\bigg)_t\,\mathrm{d}x,
 \end{split}\end{equation}
where
\begin{equation} \begin{split}\label{r4}
A&=\int_0^T \int_\Omega \psi(t)\Big(\overline{(\rho \mathbf{u})_t} - (\rho \overline{\mathbf{u}})_t\Big)\cdot\overline{\mathbf{u}}\,\mathrm{d}x\,\mathrm{d}t, \\
B&=\int_0^T \int_\Omega \psi(t) \Big(\overline{\mathrm{div}(\rho\mathbf{u}\otimes\mathbf{u})} - \mathrm{div}(\rho\mathbf{u}\otimes\overline{\mathbf{u}})\Big)\cdot\overline{\mathbf{u}}\,\mathrm{d}x \,\mathrm{d}t.
\end{split}\end{equation}

Next, we estimate the third term in \eqref{r2} as follows
\begin{equation}\begin{split}\label{r8}
&\int_0^T \int_\Omega \psi(t)\overline{\nabla P}\cdot\overline{\mathbf{u}}\,\mathrm{d}x \,\mathrm{d}t\\
%=& \int_0^T \int_\Omega \psi(t) \overline{\nabla (a \rho^\gamma)}\cdot\overline{\mathbf{u}}\,\mathrm{d}x \,\mathrm{d}t \\
=&\frac{\gamma}{\gamma - 1}\int_0^T \int_\Omega \psi(t)\overline{\rho \nabla (\rho^{\gamma-1})}\cdot\overline{\mathbf{u}}\,\mathrm{d}x \,\mathrm{d}t \\
=&\frac{\gamma}{\gamma - 1}\int_0^T \int_\Omega \psi(t)\Big(\overline{\rho\nabla(\rho^{\gamma-1})} - \rho \nabla(\rho^{\gamma-1})\Big)\cdot\overline{\mathbf{u}}\,\mathrm{d}x \,\mathrm{d}t
+ \frac{\gamma}{\gamma - 1}\int_0^T \int_\Omega \psi(t)\rho\nabla(\rho^{\gamma-1})\cdot\overline{\mathbf{u}}\,\mathrm{d}x \,\mathrm{d}t \\
%=& D +  \frac{a\gamma}{\gamma - 1}\int_0^T \int_\Omega \psi(t)\rho\overline{\nabla(\rho^{\gamma-1})}\cdot\overline{\mathbf{u}}\,\mathrm{d}x \,\mathrm{d}t \\
%=& D_1 - \frac{\gamma}{\gamma - 1}\int_0^T \int_\Omega \psi(t)\mathrm{div}(\rho \overline{\mathbf{u}}) \rho^{\gamma-1}\,\mathrm{d}x \,\mathrm{d}t\\
=& D_1 + \frac{\gamma}{\gamma - 1}\int_0^T \int_\Omega \psi(t)\Big(\overline{\mathrm{div}(\rho \mathbf{u})}-\mathrm{div}(\rho \overline{\mathbf{u}}) \Big)\rho^{\gamma-1}\,\mathrm{d}x \,\mathrm{d}t
+ \frac{\gamma}{\gamma - 1}\int_0^T \int_\Omega \psi(t)\overline{\rho_t}\rho^{\gamma-1}\,\mathrm{d}x \,\mathrm{d}t \\
=& D_1 + D_2 +\frac{\gamma}{\gamma - 1}\int_0^T \int_\Omega \psi(t)(\overline{\rho_t} - \rho_t)\rho^{\gamma-1}\,\mathrm{d}x \,\mathrm{d}t+ \frac{\gamma}{\gamma - 1}\int_0^T \int_\Omega \psi(t)\rho_t\rho^{\gamma-1}\,\mathrm{d}x \,\mathrm{d}t\\
=& D_1+ D_2+D_3+ \frac{1}{\gamma - 1}\int_0^T \int_\Omega \psi(t)(\rho^\gamma)_t\,\mathrm{d}x \,\mathrm{d}t,
\end{split}\end{equation}
where
\begin{equation}\begin{split}
D_1&=\frac{\gamma}{\gamma - 1}\int_0^T \int_\Omega \psi(t)\Big(\overline{\rho\nabla(\rho^{\gamma-1})} - \rho \nabla(\rho^{\gamma-1})\Big)\cdot\overline{\mathbf{u}}\,\mathrm{d}x \,\mathrm{d}t,\\
D_2&=\frac{\gamma}{\gamma - 1}\int_0^T \int_\Omega \psi(t)\Big(\overline{\mathrm{div}(\rho \mathbf{u})}-\mathrm{div}(\rho \overline{\mathbf{u}}) \Big)\overline{\rho^{\gamma-1}}\,\mathrm{d}x \,\mathrm{d}t,\\
D_3&= \frac{\gamma}{\gamma - 1}\int_0^T \int_\Omega \psi(t)(\overline{\rho_t} - \rho_t)\overline{\rho^{\gamma-1}}\,\mathrm{d}x \,\mathrm{d}t.
\end{split} \end{equation}

For the fifth item and the sixth item in \eqref{r2}, we have
\begin{equation}\begin{split} \label{16}
&-\int_0^T \int_\Omega \psi(t)\overline{\mu \Delta{\mathbf{u}}}\cdot\overline{\mathbf{u}}\,\mathrm{d}x \,\mathrm{d}t
-\int_0^T \int_\Omega \psi(t)\Big(\overline{(\lambda + \mu)\nabla(\mathrm{div}\mathbf{u})}\Big)\cdot\overline{\mathbf{u}}\,\mathrm{d}x \,\mathrm{d}t  \\
%=& \int_0^T \int_\Omega \mu \psi(t) \overline{\mathrm{div}(\nabla\mathbf{u})}\cdot\overline{\mathbf{u}}\,\mathrm{d}x \,\mathrm{d}t \\
%=&-\int_0^T \int_\Omega \mu \psi(t)\mathrm{div}(\overline{\nabla\mathbf{u}})\cdot\overline{\mathbf{u}}\,\mathrm{d}x \,\mathrm{d}t
%-\int_0^T \int_\Omega (\lambda + \mu) \psi(t) \nabla(\overline{\mathrm{div}\mathbf{u}})\cdot\overline{\mathbf{u}}\,\mathrm{d}x \,\mathrm{d}t \\
=&\int_0^T \int_\Omega \mu \psi(t)\mid\overline{\nabla\mathbf{u}}\mid^2\,\mathrm{d}x \,\mathrm{d}t
+\int_0^T \int_\Omega (\lambda + \mu) \psi(t)\mid\overline{\mathrm{div}\mathbf{u}}\mid^2\,\mathrm{d}x \,\mathrm{d}t.
\end{split} \end{equation}

%Next, the sixth item in \eqref{r2} is as follws:
%\begin{equation}\begin{split}
%\int_0^T& \int_\Omega \psi(t)\Big(\overline{(\lambda + \mu)\nabla(\mathrm{div}\mathbf{u})}\Big)\cdot\overline{\mathbf{u}}\,\mathrm{d}x \,\mathrm{d}t \\
%=&\int_0^T \int_\Omega (\lambda + \mu) \psi(t) \nabla(\overline{\mathrm{div}\mathbf{u}})\cdot\overline{\mathbf{u}}\,\mathrm{d}x \,\mathrm{d}t \\
%=&-\int_0^T \int_\Omega (\lambda + \mu) \psi(t)\mid\overline{\mathrm{div}\mathbf{u}}\mid^2\,\mathrm{d}x \,\mathrm{d}t.
%\end{split} \end{equation}
Finally, we handle the fourth item in \eqref{r2}.
\begin{equation}\begin{split}\label{17}
&-\int_0^T \int_\Omega \psi(t) \Big(\overline{(\nabla \times \mathbf{H}) \times \mathbf{H}}\Big)\cdot\overline{\mathbf{u}}\,\mathrm{d}x \,\mathrm{d}t \\
%=& \int_0^T \int_\Omega \psi(t)\Bigg(\overline{(\mathbf{H}\cdot\nabla)\mathbf{H} - \frac12\nabla|\mathbf{H}|^2}\Bigg)\cdot\overline{\mathbf{u}}\,\mathrm{d}x \,\mathrm{d}t \\
=& -\int_0^T \int_\Omega \psi(t)\overline{(\mathbf{H}\cdot\nabla)\mathbf{H}}\cdot\overline{\mathbf{u}}\,\mathrm{d}x \,\mathrm{d}t + \frac12\int_0^T \int_\Omega \psi(t)\overline{\nabla\mid\mathbf{H}\mid^2}\cdot\overline{\mathbf{u}}\,\mathrm{d}x \,\mathrm{d}t \\
%=& -\int_0^T \int_\Omega\psi(t)\Big(\overline{(\mathbf{H}\cdot\nabla)\mathbf{H} + \mathrm{div}(\mathbf{H})\mathbf{H}}\Big)\cdot\overline{\mathbf{u}}\,\mathrm{d}x \,\mathrm{d}t + \frac12\int_0^T \int_\Omega \psi(t)\overline{\nabla\mid\mathbf{H}\mid^2}\cdot\overline{\mathbf{u}}\,\mathrm{d}x \,\mathrm{d}t\\
=& -\int_0^T \int_\Omega\psi(t)\overline{\mathrm{div}(\mathbf{H}\otimes \mathbf{H})}\cdot\overline{\mathbf{u}}\,\mathrm{d}x \,\mathrm{d}t
- \frac12\int_0^T \int_\Omega \psi(t)\mathrm{div}\overline{\mathbf{u}}\overline{\mathbf{H}\cdot\mathbf{H}}\,\mathrm{d}x \,\mathrm{d}t.
\end{split}\end{equation}
The first term in the last equality of \eqref{17} shows that
\begin{equation}\label{3.8}\begin{split}
&-\int_0^T \int_\Omega\psi(t)\overline{\mathrm{div}(\mathbf{H}\otimes \mathbf{H})}\cdot\overline{\mathbf{u}}\,\mathrm{d}x \,\mathrm{d}t\\
=&-\int_0^T \int_\Omega\psi(t)\Big(\overline{\mathrm{div}(\mathbf{H}\otimes \mathbf{H})}-\mathrm{div}(\mathbf{H}\otimes\overline{\mathbf{H}})\Big)\cdot\overline{\mathbf{u}}\,\mathrm{d}x \,\mathrm{d}t-  \int_0^T \int_\Omega\psi(t)\mathrm{div}(\mathbf{H}\otimes\overline{\mathbf{H}})\cdot\overline{\mathbf{u}}\,\mathrm{d}x \,\mathrm{d}t\\
%=&G_1 - \int_0^T \int_\Omega\psi(t)\mathrm{div}(\mathbf{H}\otimes\overline{\mathbf{H}})\cdot\overline{\mathbf{u}}\,\mathrm{d}x \,\mathrm{d}t\\
%=&G_1 - \int_0^T \int_\Omega\psi(t)\Big((\mathbf{H}\cdot\nabla)\overline{\mathbf{H}} + \mathrm{div}(\mathbf{H})\overline{\mathbf{H}}\Big)\cdot\overline{\mathbf{u}}\,\mathrm{d}x\,\mathrm{d}t\\
=&I_1 - \int_0^T \int_\Omega\psi(t)\Big((\mathbf{H}\cdot\nabla)\overline{\mathbf{H}}\Big)\cdot\overline{\mathbf{u}}\,\mathrm{d}x\,\mathrm{d}t .
\end{split}\end{equation}
And the second term in the last equality of \eqref{17} shows that
\begin{equation}\label{3.9}\begin{split}
%\frac12&\int_0^T \int_\Omega \psi(t)\overline{\nabla\mid\mathbf{H}\mid^2}\cdot\overline{\mathbf{u}}\,\mathrm{d}x \,\mathrm{d}t\\
&-\frac12\int_0^T \int_\Omega \psi(t)\mathrm{div}\overline{\mathbf{u}}\overline{\mathbf{H}\cdot\mathbf{H}}\,\mathrm{d}x \,\mathrm{d}t\\
=&-\frac12\int_0^T \int_\Omega \psi(t)\mathrm{div}\overline{\mathbf{u}}(\overline{\mathbf{H}\cdot\mathbf{H}} - \mathbf{H}\cdot \overline{\mathbf{H}})\,\mathrm{d}x \,\mathrm{d}t - \frac12\int_0^T \int_\Omega \psi(t)\mathrm{div}\mathbf{\overline{u}}(\mathbf{H}\cdot\overline{\mathbf{H}})\,\mathrm{d}x\,\mathrm{d}t\\
=&\frac12\int_0^T \int_\Omega \psi(t)\overline{\mathbf{u}}\cdot(\overline{\nabla(\mathbf{H}\cdot\mathbf{H})} - \nabla(\mathbf{H}\cdot \overline{\mathbf{H}}))\,\mathrm{d}x \,\mathrm{d}t- \frac12\int_0^T \int_\Omega \psi(t)\mathrm{div}\mathbf{\overline{u}}(\mathbf{H}\cdot\overline{\mathbf{H}})\,\mathrm{d}x\,\mathrm{d}t\\
=&I_2- \frac12\int_0^T \int_\Omega \psi(t)\mathrm{div}\mathbf{\overline{u}}(\mathbf{H}\cdot\overline{\mathbf{H}})\,\mathrm{d}x\,\mathrm{d}t.
\end{split}\end{equation}

Substituting (\ref{3.8}) and (\ref{3.9}) into \eqref{17}, we obtain
\begin{equation}\begin{split}\label{18}
&-\int_0^T\int_\Omega \psi(t) \Big(\overline{(\nabla \times \mathbf{H}) \times \mathbf{H}}\Big)\cdot\overline{\mathbf{u}}\,\mathrm{d}x \,\mathrm{d}t \\
=&I_1 + I_2 - \int_0^T \int_\Omega\psi(t)\Big((\mathbf{H}\cdot\nabla)\overline{\mathbf{H}}\Big)\cdot\overline{\mathbf{u}}\,\mathrm{d}x\,\mathrm{d}t - \frac12\int_0^T \int_\Omega \psi(t)\mathrm{div}\mathbf{\overline{u}}(\mathbf{H}\cdot\overline{\mathbf{H}})\,\mathrm{d}x\,\mathrm{d}t,
\end{split}\end{equation}
where
\begin{equation}\begin{split}
I_1 &=  \int_0^T \int_\Omega\psi(t)\Big(\mathrm{div}(\mathbf{H}\otimes\overline{\mathbf{H}})
-\overline{\mathrm{div}(\mathbf{H}\otimes\mathbf{H})}\Big)\cdot\overline{\mathbf{u}}\,\mathrm{d}x \,\mathrm{d}t,\\
I_2 &= \frac12\int_0^T \int_\Omega \psi(t)\overline{\mathbf{u}}\cdot\Big(\overline{\nabla(\mathbf{H}\cdot\mathbf{H})} - \nabla(\mathbf{H}\cdot \overline{\mathbf{H}})\Big)\,\mathrm{d}x \,\mathrm{d}t.
\end{split} \end{equation}

%So now we can get the estimation of \eqref{r2} as follows:
Combining \eqref{r3}, \eqref{r8}, \eqref{16} with \eqref{18}, we can get the equality of \eqref{r2} as follows
\begin{equation}\begin{split}\label{r23}
%&\int_0^T \int_\Omega \psi(t)\overline{\mathbf{u}}\cdot \Big(\overline{ (\rho\mathbf{u})_t + \mathrm{div}(\rho\mathbf{u} \otimes \mathbf{u}) + \nabla{P} - (\nabla \times \mathbf{H}) \times \mathbf{H} - \mu \Delta{\mathbf{u}} - (\lambda + \mu)\nabla(\mathrm{div}\mathbf{u})}\Big)\,\mathrm{d}x\,\mathrm{d}t\\
&\int_0^T \int_\Omega \psi(t)\bigg(\frac{1}{2} \rho\mid\overline{\mathbf{u}}\mid^2 +\frac{\rho^\gamma}{\gamma - 1} \bigg)_t\,\mathrm{d}x\,\mathrm{d}t + \int_0^T \int_\Omega \psi(t)\Big(\mu \mid\overline{\nabla\mathbf{u}}\mid^2+(\lambda + \mu)\mid\overline{\mathrm{div}\mathbf{u}}\mid^2\Big)\,\mathrm{d}x \,\mathrm{d}t \\
&- \int_0^T \int_\Omega\psi(t)\Big((\mathbf{H}\cdot\nabla)\overline{\mathbf{H}}\Big)\cdot\overline{\mathbf{u}}\,\mathrm{d}x\,\mathrm{d}t - \frac12\int_0^T \int_\Omega \psi(t)\mathrm{div}\mathbf{\overline{u}}(\mathbf{H}\cdot\overline{\mathbf{H}})\,\mathrm{d}x\,\mathrm{d}t+\mathbf{R}_{\epsilon}+I_1+I_2=0,
\end{split}\end{equation}
where
\begin{equation}
\mathbf{R}_{\epsilon}=A + B + D_1+ D_2 + D_3.
\end{equation}

%So from the equation \eqref{r2} we have the following:
%\begin{equation}\begin{split}\label{r23}
%&\int_0^T \int_\Omega \psi(t)\bigg(\frac{1}{2} \rho\mid\overline{\mathbf{u}}\mid^2 +\frac{a}{\gamma - 1} \rho^\gamma\bigg)_t\,\mathrm{d}x\,\mathrm{d}t + % \int_0^T \int_\Omega \mu \psi(t)\mid\overline{\nabla\mathbf{u}}\mid^2\,\mathrm{d}x \,\mathrm{d}t + \int_0^T \int_\Omega (\lambda + \mu) \psi(t)\mid\overline{\mathrm{div}\mathbf{u}}\mid^2\,\mathrm{d}x \,\mathrm{d}t \\
%&- \int_0^T \int_\Omega\psi(t)\Big((\mathbf{H}\cdot\nabla)\overline{\mathbf{H}}\Big)\cdot\overline{\mathbf{u}}\,\mathrm{d}x\,\mathrm{d}t - \frac12\int_0^T \int_\Omega \psi(t)\mathrm{div}(\mathbf{\overline{u}})(\mathbf{H}\cdot\overline{\mathbf{H}})\,\mathrm{d}x\,\mathrm{d}t +\mathbf{R}_{\epsilon}- G_{1} + G_{2}= 0.
%\end{split}\end{equation}

Now we are in a position to handle $\eqref{2dbu-E1.1}_3$. Here we introduce a new function $\Theta = \overline{\psi(t)\overline{\mathbf{H}}}$ as a test function of $\eqref{2dbu-E1.1}_3$. Then we get
\begin{equation}\begin{split}\label{r24}
 \int_0^T \int_\Omega\psi(t)\Big(\overline{\mathbf{H}_t + \nabla\times(\nu\nabla\times\mathbf{H}) - \nabla\times(\mathbf{u}\times\mathbf{H})}\Big)\cdot\overline{\mathbf{H}}\,\mathrm{d}x \,\mathrm{d}t = 0.
\end{split} \end{equation}

The first term in \eqref{r24} shows that
\begin{equation}\begin{split}
\int_0^T \int_\Omega\psi(t)\overline{\mathbf{H}_t}\cdot\overline{\mathbf{H}}\,\mathrm{d}x \,\mathrm{d}t = \frac12\int_0^T \int_\Omega\psi(t)
(\mid\overline{\mathbf{H}}\mid^2)_t\,\mathrm{d}x \,\mathrm{d}t.
\end{split} \end{equation}

Similarly, the second term in \eqref{r24} yields
\begin{equation}\begin{split}
\int_0^T \int_\Omega\psi(t)\overline{\nabla\times(\nu\nabla\times\mathbf{H})}\cdot\overline{\mathbf{H}}\,\mathrm{d}x \,\mathrm{d}t
%=& \nu\int_0^T \int_\Omega\psi(t)\overline{\nabla(\mathrm{div}\mathbf{H}) - \Delta\mathbf{H}}\cdot\overline{\mathbf{H}}\,\mathrm{d}x \,\mathrm{d}t\\
=& -\nu\int_0^T \int_\Omega\psi(t)\overline{\Delta\mathbf{H}}\cdot\overline{\mathbf{H}}\,\mathrm{d}x \,\mathrm{d}t\\
%=& -\nu\int_0^T \int_\Omega\psi(t)\overline{\mathrm{div}(\nabla\mathbf{H})}\cdot\overline{\mathbf{H}}\,\mathrm{d}x \,\mathrm{d}t\\
%=&-\nu\int_0^T \int_\Omega\psi(t)\mathrm{div}(\overline{\nabla\mathbf{H}})\cdot\overline{\mathbf{H}}\,\mathrm{d}x \,\mathrm{d}t\\
=&\nu\int_0^T \int_\Omega\psi(t)\mid\overline{\nabla\times\mathbf{H}}\mid^2\,\mathrm{d}x \,\mathrm{d}t.
\end{split} \end{equation}

Finally, the last term in \eqref{r24} shows that
\begin{equation}\begin{split}\label{28}
& -\int_0^T\int_\Omega\psi(t)\overline{\mathbf{H}}\cdot\overline{\nabla\times(\mathbf{u}\times\mathbf{H})}\,\mathrm{d}x \,\mathrm{d}t\\
=&-\int_0^T \int_\Omega\psi(t)\overline{\mathbf{H}}\cdot\overline{\mathbf{u}(\mathrm{div}\mathbf{H}) - \mathbf{H}(\mathrm{div}\mathbf{u}) + (\mathbf{H}\cdot\nabla)\mathbf{u} - (\mathbf{u}\cdot\nabla)\mathbf{H}}\,\mathrm{d}x \,\mathrm{d}t\\
=&-\int_0^T \int_\Omega\psi(t)\overline{\mathbf{H}}\cdot\overline{(\mathbf{H}\cdot\nabla)\mathbf{u} - \mathbf{H}(\mathrm{div}\mathbf{u})  - (\mathbf{u}\cdot\nabla)\mathbf{H}}\,\mathrm{d}x \,\mathrm{d}t\\
=&-\int_0^T \int_\Omega\psi(t)\overline{\mathbf{H}}\cdot\overline{(\mathbf{H}\cdot\nabla)\mathbf{u}}\,\mathrm{d}x \,\mathrm{d}t + \int_0^T \int_\Omega\psi(t)\overline{\mathbf{H}}\cdot\overline{ (\mathbf{u}\cdot\nabla)\mathbf{H} +(\mathrm{div}\mathbf{u})\mathbf{H} }\,\mathrm{d}x \,\mathrm{d}t\\
%=&\int_0^T \int_\Omega\psi(t)\overline{\mathbf{H}}\cdot\overline{(\mathbf{H}\cdot\nabla)\mathbf{u}+(\mathrm{div}\mathbf{H})\mathbf{u}}\,\mathrm{d}x \,\mathrm{d}t- \int_0^T \int_\Omega\psi(t)\overline{\mathbf{H}}\cdot\overline{ (\mathbf{u}\cdot\nabla)\mathbf{H} +(\mathrm{div}\mathbf{u})\mathbf{H} }\,\mathrm{d}x \,\mathrm{d}t\\
=&\D-\int_0^T \int_\Omega\psi(t)\overline{\mathbf{H}}\cdot\overline{\mathrm{div}(\mathbf{H}\otimes\mathbf{u})}\,\mathrm{d}x \,\mathrm{d}t + \int_0^T \int_\Omega\psi(t)\overline{\mathbf{H}}\cdot\overline{\mathrm{div}(\mathbf{u}\otimes\mathbf{H})}\,\mathrm{d}x \,\mathrm{d}t. \end{split} \end{equation}
The first term in the last equality of \eqref{28} shows that
\begin{equation}\begin{split}
&-\int_0^T \int_\Omega\psi(t)\overline{\mathbf{H}}\cdot\overline{\mathrm{div}(\mathbf{H}\otimes\mathbf{u})}\,\mathrm{d}x \,\mathrm{d}t\\
=&-\int_0^T\int_\Omega\psi(t)\overline{\mathbf{H}}\cdot\Big(\overline{\mathrm{div}(\mathbf{H}\otimes\mathbf{u})}-\mathrm{div}
(\mathbf{H}\otimes\overline{\mathbf{u}})\Big)\,\mathrm{d}x \,\mathrm{d}t - \int_0^T \int_\Omega\psi(t)\overline{\mathbf{H}}\cdot \mathrm{div}(\mathbf{H}\otimes\overline{\mathbf{u}})\,\mathrm{d}x \,\mathrm{d}t\\
%=&I_3 -\int_0^T \int_\Omega\psi(t)\overline{\mathbf{H}}\cdot \mathrm{div}(\mathbf{H}\otimes\overline{\mathbf{u}})\,\mathrm{d}x \,\mathrm{d}t\\
%=&I_1  + \int_0^T \int_\Omega\psi(t)\overline{\mathbf{H}}\cdot\Big( (\mathbf{H}\cdot\nabla)\mathbf{\overline{u}} + (\mathrm{div}\mathbf{H})\overline{\mathbf{u}}\Big)\,\mathrm{d}x \,\mathrm{d}t \\
=&I_3 -  \int_0^T \int_\Omega\psi(t)\overline{\mathbf{H}}\cdot\Big((\mathbf{H}\cdot\nabla)\mathbf{\overline{u}}\Big)\,\mathrm{d}x \,\mathrm{d}t.\\
\end{split} \end{equation}
And the second term in the last equality of \eqref{28} shows that
\begin{equation}\begin{split}
&\int_0^T \int_\Omega\psi(t)\overline{\mathbf{H}}\cdot\overline{\mathrm{div}(\mathbf{u}\otimes\mathbf{H})}\,\mathrm{d}x \,\mathrm{d}t\\
=&\int_0^T \int_\Omega\psi(t)\overline{\mathbf{H}}\cdot\Big(\overline{\mathrm{div}(\mathbf{u}\otimes\mathbf{H})}-\mathrm{div}(\overline{\mathbf{u}}\otimes\mathbf{H})\Big)\,\mathrm{d}x \,\mathrm{d}t + \int_0^T \int_\Omega\psi(t)\overline{\mathbf{H}}\cdot\mathrm{div}(\overline{\mathbf{u}}\otimes\mathbf{H})\,\mathrm{d}x \,\mathrm{d}t\\
=&I_4 +  \int_0^T \int_\Omega\psi(t)\overline{\mathbf{H}}\cdot\mathrm{div}(\overline{\mathbf{u}}\otimes\mathbf{H})\,\mathrm{d}x \,\mathrm{d}t\\
%=&I_2 + \int_0^T \int_\Omega\psi(t)\overline{\mathbf{H}}\cdot\Big((\overline{\mathbf{u}}\cdot\nabla)\mathbf{H} + (\mathrm{div}\overline{\mathbf{u}})\mathbf{H} \Big)\,\mathrm{d}x \,\mathrm{d}t\\
=&I_4 +  \int_0^T \int_\Omega\psi(t)\overline{\mathbf{H}}\cdot\Big((\overline{\mathbf{u}}\cdot\nabla)\mathbf{H}\Big)\,\mathrm{d}x \,\mathrm{d}t +\int_0^T \int_\Omega\psi(t)\mathrm{div}\overline{\mathbf{u}}(\mathbf{H}\cdot\overline{\mathbf{H}})\,\mathrm{d}x \,\mathrm{d}t.\\
\end{split} \end{equation}
 Substituting the above two equalities into \eqref{28}, we have
\begin{equation}\begin{split}
&-\int_0^T\int_\Omega\psi(t)\overline{\mathbf{H}}\cdot\overline{\nabla\times(\mathbf{u}\times\mathbf{H})}\,\mathrm{d}x \,\mathrm{d}t\\
=&I_3 + I_4 - \int_0^T \int_\Omega\psi(t)\overline{\mathbf{H}}\cdot\Big((\mathbf{H}\cdot\nabla)\mathbf{\overline{u}}\Big)\,\mathrm{d}x \,\mathrm{d}t + \int_0^T \int_\Omega\psi(t)\overline{\mathbf{H}}\cdot\Big((\overline{\mathbf{u}}\cdot\nabla)\mathbf{H}\Big)\,\mathrm{d}x \,\mathrm{d}t \\
&+\int_0^T \int_\Omega\psi(t)\mathrm{div}\overline{\mathbf{u}}(\mathbf{H}\cdot\overline{\mathbf{H}})\,\mathrm{d}x \,\mathrm{d}t.
\end{split} \end{equation}

Recalling \eqref{r24}, we have
\begin{equation}\begin{split}\label{r28}
&\frac12\int_0^T \int_\Omega\psi(t)(\mid\overline{\mathbf{H}}\mid^2)_t\,\mathrm{d}x \,\mathrm{d}t +\nu\int_0^T \int_\Omega\psi(t)\mid\overline{\nabla\times\mathbf{H}}\mid^2\,\mathrm{d}x \,\mathrm{d}t
+ \int_0^T \int_\Omega\psi(t)\overline{\mathbf{H}}\cdot\Big((\overline{\mathbf{u}}\cdot\nabla)\mathbf{H}\Big)\,\mathrm{d}x \,\mathrm{d}t\\
 &+ \int_0^T \int_\Omega\psi(t)\mathrm{div}\overline{\mathbf{u}}(\mathbf{H}\cdot\overline{\mathbf{H}})\,\mathrm{d}x \,\mathrm{d}t
-  \int_0^T \int_\Omega\psi(t)\overline{\mathbf{H}}\cdot\Big((\mathbf{H}\cdot\nabla)\mathbf{\overline{u}}\Big)\,\mathrm{d}x \,\mathrm{d}t+I_3 + I_4= 0,
\end{split} \end{equation}
where
\begin{equation}\begin{split}
I_3 &= \int_0^T \int_\Omega\psi(t)\overline{\mathbf{H}}\cdot\Big(\mathrm{div}(\mathbf{H}\otimes\overline{\mathbf{u}})
-\overline{\mathrm{div}(\mathbf{H}\otimes\mathbf{u})}\Big)\,\mathrm{d}x \,\mathrm{d}t,\\
I_4 &=  \int_0^T \int_\Omega\psi(t)\overline{\mathbf{H}}\cdot\Big(\overline{\mathrm{div}(\mathbf{u}\otimes\mathbf{H})}
-\mathrm{div}(\overline{\mathbf{u}}\otimes\mathbf{H})\Big)\,\mathrm{d}x \,\mathrm{d}t.
\end{split} \end{equation}

Combining \eqref{r23} with \eqref{r28}, we have
\begin{equation}\begin{split}\label{r30}
&\int_0^T \int_\Omega \psi(t)\bigg(\frac{1}{2} \rho\mid\overline{\mathbf{u}}\mid^2 +\frac{\rho^\gamma}{\gamma - 1}  +\frac{1}{2} \mid\overline{\mathbf{H}}\mid^2\bigg)_t\,\mathrm{d}x\,\mathrm{d}t +\nu\int_0^T \int_\Omega\psi(t)\mid\overline{\nabla\times\mathbf{H}}\mid^2\,\mathrm{d}x \,\mathrm{d}t \\
&+\int_0^T \int_\Omega \mu \psi(t)\mid\overline{\nabla\mathbf{u}}\mid^2\,\mathrm{d}x \,\mathrm{d}t + \int_0^T \int_\Omega (\lambda + \mu) \psi(t)\mid\overline{\mathrm{div}\mathbf{u}}\mid^2\,\mathrm{d}x \,\mathrm{d}t\\
&- \int_0^T \int_\Omega\psi(t)\Big((\mathbf{H}\cdot\nabla)\overline{\mathbf{H}}\Big)\cdot\overline{\mathbf{u}}\,\mathrm{d}x\,\mathrm{d}t- \int_0^T \int_\Omega\psi(t)\overline{\mathbf{H}}\cdot\Big((\mathbf{H}\cdot\nabla)\mathbf{\overline{u}}\Big)\,\mathrm{d}x \,\mathrm{d}t\\
&+ \frac12\int_0^T \int_\Omega \psi(t)\mathrm{div}\mathbf{\overline{u}}(\mathbf{H}\cdot\overline{\mathbf{H}})\,\mathrm{d}x\,\mathrm{d}t  + \int_0^T \int_\Omega\psi(t)\overline{\mathbf{H}}\cdot\Big((\overline{\mathbf{u}}\cdot\nabla)\mathbf{H}\Big)\,\mathrm{d}x \,\mathrm{d}t +R_{\epsilon} +I_{\epsilon}=0,
\end{split} \end{equation}
where
\begin{equation}
I_{\epsilon} = I_1+I_2+I_3+I_4.
\end{equation}

In equation \eqref{r30}, we continue to estimate the last four terms as follows\\
On the one hand, we have
\begin{equation}\begin{split}
 - &\int_0^T \int_\Omega\psi(t)\Big((\mathbf{H}\cdot\nabla)\overline{\mathbf{H}}\Big)\cdot\overline{\mathbf{u}}\,\mathrm{d}x\,\mathrm{d}t-  \int_0^T \int_\Omega\psi(t)\overline{\mathbf{H}}\cdot\Big((\mathbf{H}\cdot\nabla)\mathbf{\overline{u}}\Big)\,\mathrm{d}x \,\mathrm{d}t
%=&\int_0^T \int_\Omega\psi(t)\mathrm{div}(\mathbf{H})(\overline{\mathbf{H}}\cdot\overline{\mathbf{u}})\,\mathrm{d}x \,\mathrm{d}t + \int_0^T \int_\Omega\psi(t)\overline{\mathbf{H}}\cdot\Big((\mathbf{H}\cdot\nabla)\mathbf{\overline{u}}\Big)\,\mathrm{d}x \,\mathrm{d}t - \int_0^T \int_\Omega\psi(t)\overline{\mathbf{H}}\cdot\Big((\mathbf{H}\cdot\nabla)\mathbf{\overline{u}}\Big)\,\mathrm{d}x \,\mathrm{d}t\\
=0.
\end{split} \end{equation}
On the other hand, we deduce
\begin{align}
%& -\frac12\int_0^T \int_\Omega \psi(t)\mathrm{div}(\mathbf{\overline{u}})\Big(\mathbf{H}\cdot\overline{\mathbf{H}}\Big)\,\mathrm{d}x\,\mathrm{d}t
%+ \int_0^T \int_\Omega\psi(t)\overline{\mathbf{H}}\cdot\Big((\overline{\mathbf{u}}\cdot\nabla)\mathbf{H}\Big)\,\mathrm{d}x \,\mathrm{d}t +  \int_0^T \int_\Omega\psi(t)\overline{\mathbf{H}}\cdot\big[(\mathrm{div}\overline{\mathbf{u}})\mathbf{H} \big]\,\mathrm{d}x \,\mathrm{d}t\\
&\int_0^T \int_\Omega\psi(t)\overline{\mathbf{H}}\cdot\Big((\overline{\mathbf{u}}\cdot\nabla)\mathbf{H}\Big)\,\mathrm{d}x \,\mathrm{d}t+ \frac12\int_0^T \int_\Omega \psi(t)\mathrm{div}\mathbf{\overline{u}}(\mathbf{H}\cdot\overline{\mathbf{H}})\,\mathrm{d}x\,\mathrm{d}t\notag\\
=&\int_0^T \int_\Omega\psi(t)\overline{\mathbf{H}}\cdot\Big((\overline{\mathbf{u}}\cdot\nabla)\mathbf{H}\Big)\,\mathrm{d}x \,\mathrm{d}t -  \frac12\int_0^T \int_\Omega \psi(t)\overline{\mathbf{u}}\cdot\nabla(\mathbf{H}\cdot\overline{\mathbf{H}})\,\mathrm{d}x\,\mathrm{d}t\notag\\
%=&\int_0^T \int_\Omega\psi(t)\overline{\mathbf{H}}\cdot\Big((\overline{\mathbf{u}}\cdot\nabla)\mathbf{H}\Big)\,\mathrm{d}x \,\mathrm{d}t -  \frac12\int_0^T \int_\Omega \psi(t)\Big((\overline{\mathbf{u}}\cdot\nabla)\mathbf{H}\Big)\cdot\overline{\mathbf{H}}\,\mathrm{d}x\,\mathrm{d}t -  \frac12\int_0^T \int_\Omega \psi(t)\Big((\overline{\mathbf{u}}\cdot\nabla)\overline{\mathbf{H}}\Big)\cdot\mathbf{H}\,\mathrm{d}x\,\mathrm{d}t\\
=&\frac12\int_0^T \int_\Omega\psi(t)\overline{\mathbf{H}}\cdot\Big((\overline{\mathbf{u}}\cdot\nabla)\mathbf{H}\Big)\,\mathrm{d}x \,\mathrm{d}t
-  \frac12\int_0^T \int_\Omega \psi(t)\mathbf{H}\cdot\Big((\overline{\mathbf{u}}\cdot\nabla)\overline{\mathbf{H}}\Big)\,\mathrm{d}x\,\mathrm{d}t\\
=&\frac12\int_0^T \int_\Omega\psi(t)\overline{\mathbf{H}}\cdot\Big((\overline{\mathbf{u}}\cdot\nabla)\mathbf{H}\Big)\,\mathrm{d}x \,\mathrm{d}t - \frac12\int_0^T \int_\Omega\psi(t)\mathbf{H}\cdot\Big((\overline{\mathbf{u}}\cdot\nabla)\mathbf{H}\Big)\,\mathrm{d}x \,\mathrm{d}t\notag\\
&+\frac12\int_0^T \int_\Omega\psi(t)\mathbf{H}\cdot\Big((\overline{\mathbf{u}}\cdot\nabla)\mathbf{H}\Big)\,\mathrm{d}x \,\mathrm{d}t -  \frac12\int_0^T \int_\Omega \psi(t)\mathbf{H}\cdot\Big((\overline{\mathbf{u}}\cdot\nabla)\overline{\mathbf{H}}\Big)\,\mathrm{d}x\,\mathrm{d}t\notag\\
%=&\frac12\int_0^T \int_\Omega\psi(t)(\overline{\mathbf{H}} - \mathbf{H})\cdot\Big((\overline{\mathbf{u}}\cdot\nabla)\mathbf{H}\Big)\,\mathrm{d}x \,\mathrm{d}t  + \frac12\int_0^T \int_\Omega\psi(t)\mathbf{H}\cdot\Big((\overline{\mathbf{u}}\cdot\nabla)(\overline{\mathbf{H}}-\mathbf{H})\Big)\,\mathrm{d}x \,\mathrm{d}t\\
=&J_1 + J_2,\notag
\end{align}
where
\begin{equation}\begin{split}
J_1 &= \frac12\int_0^T \int_\Omega\psi(t)(\overline{\mathbf{H}} - \mathbf{H})\cdot\Big((\overline{\mathbf{u}}\cdot\nabla)\mathbf{H}\Big)\,\mathrm{d}x \,\mathrm{d}t,\\
J_2 &= \frac12\int_0^T \int_\Omega\psi(t)\mathbf{H}\cdot\Big((\overline{\mathbf{u}}\cdot\nabla)(\mathbf{H}-\overline{\mathbf{H}})\Big)\,\mathrm{d}x \,\mathrm{d}t.
\end{split} \end{equation}

By the above equalities and integration by parts, it yields
\begin{equation}\begin{split}\label{r40}
-&\int_0^T \int_\Omega \psi_t\bigg(\frac{1}{2} \rho\mid\overline{\mathbf{u}}\mid^2 +\frac{\rho^\gamma}{\gamma - 1}  + \frac{1}{2}\mid\overline{\mathbf{H}}\mid^2\bigg)\,\mathrm{d}x\,\mathrm{d}t\\
+&\int_0^T \int_\Omega\psi(t)\Big(\nu\mid\overline{\nabla\times\mathbf{H}}\mid^2+ \mu\mid\overline{\nabla\mathbf{u}}\mid^2+ (\lambda + \mu)\mid\overline{\mathrm{div}\mathbf{u}}\mid^2\Big)\,\mathrm{d}x \,\mathrm{d}t+R_{\epsilon} + \mathit{I}_{\epsilon} +J_\epsilon=0,
\end{split}\end{equation}
where $J_\epsilon=J_1+J_2.$\\

\textbf{Step 2.} Passing to the limit in \eqref{r40} as $\epsilon$ tends to zero.

Using Definition \ref{defn}, \eqref{rho} and \eqref{u}, one obtains
\begin{equation}\begin{split}
&\int_0^T \int_\Omega \frac{1}{2} \rho\mid\overline{\mathbf{u}}\mid^2\psi_t\,\mathrm{d}x \,\mathrm{d}t \to \int_0^T \int_\Omega \frac{1}{2} \rho\mid\mathbf{u}\mid^2\psi_t\,\mathrm{d}x \,\mathrm{d}t,\\
&\int_0^T \int_\Omega\frac{1}{2}\mid\overline{\mathbf{H}}\mid^2\psi_t\,\mathrm{d}x \,\mathrm{d}t \to \int_0^T \int_\Omega\frac{1}{2}\mid\mathbf{H}\mid^2\psi_t\,\mathrm{d}x \,\mathrm{d}t,\\
&\int_0^T \int_\Omega\psi(t)\nu\mid\overline{\nabla\times\mathbf{H}}\mid^2\,\mathrm{d}x \,\mathrm{d}t \to \int_0^T \int_\Omega\psi(t)\nu\mid\nabla\times\mathbf{H}\mid^2\,\mathrm{d}x \,\mathrm{d}t,\\
&\int_0^T \int_\Omega\psi(t)\mu\mid\overline{\nabla\mathbf{u}}\mid^2\,\mathrm{d}x \,\mathrm{d}t\to \int_0^T \int_\Omega\psi(t)\mu\mid\nabla\mathbf{u}\mid^2\,\mathrm{d}x \,\mathrm{d}t,\\
&\int_0^T \int_\Omega\psi(t) (\lambda + \mu)\mid\overline{\mathrm{div}\mathbf{u}}\mid^2\,\mathrm{d}x \,\mathrm{d}t \to \int_0^T \int_\Omega\psi(t) (\lambda + \mu) \mid\mathrm{div}\mathbf{u}\mid^2\,\mathrm{d}x \,\mathrm{d}t,
\end{split}\end{equation}
as $\epsilon \to 0$.\\

The next goal is to make use of Lemma \ref{limit} to prove
\begin{equation}
R_{\epsilon} + \mathit{I}_{\epsilon} +J_\epsilon\to 0, ~as ~\epsilon \to 0.
\end{equation}

Firstly, we prove $R_{\epsilon} \to 0$, as $\epsilon \to 0$. We assume that $\mathbf{u}$ is bounded in $L^p(0,T;L^q(\Omega))$.
On the one hand, due to \eqref{regu}, \eqref{rho}, we have
\begin{equation}\label{rho_t}
  \rho_{t}=-(\rho \mathrm{div}\mathbf{u}+2 \sqrt{\rho} \mathbf{u} \cdot \nabla \sqrt{\rho})\in L^{2}\left(0,T; L^{2}(\Omega)\right)+L^{p}\left(0,T ; L^{\frac{2q}{q+2}}(\Omega)\right).
\end{equation}
Thus, in view of Lemma \ref{limit}, we have
\begin{equation}\begin{split}
A+B\leq&\mid A\mid+\mid B\mid=\Big|\int_0^T \int_\Omega \psi(t)\Big(\overline{(\rho \mathbf{u})_t} - (\rho \overline{\mathbf{u}})_t\Big)\cdot\overline{\mathbf{u}}\,\mathrm{d}x\,\mathrm{d}t\Big|\\
&+\Big|\int_0^T \int_\Omega \psi(t) \Big(\overline{\mathrm{div}(\rho\mathbf{u}\otimes\mathbf{u})} - \mathrm{div}(\rho\mathbf{u}\otimes\overline{\mathbf{u}})\Big)\cdot\overline{\mathbf{u}}\,\mathrm{d}x \,\mathrm{d}t\Big|\\
%\leq&\parallel\psi(t)\parallel_{L^{\infty}(0,T)}\int_0^T \int_\Omega\mid\Big(\overline{(\rho \mathbf{u})_t} - (\rho \overline{\mathbf{u}})_t\Big)\cdot\overline{\mathbf{u}}\mid\,\mathrm{d}x\,\mathrm{d}t\\
\leq&\parallel\psi(t)\parallel_{L^{\infty}(0,T)}\int_0^T\parallel\mathbf{u}\parallel_{L^q(\Omega)}\bigg(\parallel\overline{(\rho \mathbf{u})_t} - (\rho \overline{\mathbf{u}})_t\parallel_{L^{\frac{q}{q-1}}(\Omega)}\\
&+\parallel\overline{\mathrm{div}(\rho\mathbf{u}\otimes\mathbf{u})} - \mathrm{div}(\rho\mathbf{u}\otimes\overline{\mathbf{u}})\parallel_{L^{\frac{q}{q-1}}(\Omega)}\bigg)\mathrm{d}t\\
\leq&\parallel\psi(t)\parallel_{L^{\infty}(0,T)}\int_0^T\parallel\mathbf{u}\parallel_{L^q(\Omega)}\parallel\rho_t\parallel_{L^{\frac{2q}{q+2}}(\Omega)}
\parallel\mathbf{u}\parallel_{L^{\frac{2q}{q-4}}(\Omega)}\mathrm{d}t\\
\leq&C\int_0^T\parallel\mathbf{u}\parallel^4_{L^q(\Omega)}
+\parallel\rho_t\parallel^2_{L^{\frac{2q}{q+2}}(\Omega)}
\mathrm{d}t\\
\leq&C\parallel\mathbf{u}\parallel^4_{L^p(0,T;L^q(\Omega))}
+C\parallel\rho_t\parallel^2_{L^2(0,T;L^{\frac{2q}{q+2}}(\Omega))},
\end{split}\end{equation}
for any $p\geq4$ and $q\geq6$.

Thanks to Lemma \ref{limit}, as $\epsilon \to 0$, we have
\begin{equation}\begin{split}
A+B \to 0.
\end{split}\end{equation}
For $D_2,$ we get
\begin{align}
D_2\leq\mid D_2 \mid=&\Big|\frac{\gamma}{\gamma - 1}\int_0^T \int_\Omega \psi(t)\Big(
 \overline{\mathrm{div}(\rho \mathbf{u})}-\mathrm{div}(\rho \overline{\mathbf{u}})\Big)\overline{\rho^{\gamma-1}}\,\mathrm{d}x \,\mathrm{d}t\Big|\notag\\
\leq&C\parallel\psi(t)\parallel_{L^{\infty}(0,T)}\int_0^T\parallel\overline{\rho^{\gamma-1}}\parallel_{L^\infty(\Omega)}\parallel
 \overline{\mathrm{div}(\rho \mathbf{u})}-\mathrm{div}(\rho \overline{\mathbf{u}})\parallel_{L^1(\Omega)}\mathrm{d}t\notag\\
\leq&C\parallel\psi(t)\parallel_{L^{\infty}(0,T)}\int_0^T\parallel\sqrt{\rho}\nabla\sqrt{\rho}\parallel_{L^\frac{q}{q-1}(\Omega)}\parallel
\mathbf{u}\parallel_{L^q(\Omega)}\mathrm{d}t\\
\leq&C\parallel\psi(t)\parallel_{L^{\infty}(0,T)}\int_0^T\parallel\nabla\sqrt{\rho}\parallel_{L^2(\Omega)}\parallel
\mathbf{u}\parallel_{L^q(\Omega)}\mathrm{d}t\notag\\
\leq&C\parallel\nabla\sqrt{\rho}\parallel_{L^\infty(0,T;L^2(\Omega))}\parallel
\mathbf{u}\parallel_{L^p(0,T;L^q(\Omega))}\mathrm{d}t\notag,
\end{align}
for any $q\geq6.$
Similarly, by Lemma \ref{limit}, we have that $D_2$ converges to zero, as $\epsilon\rightarrow0$.

For $D_1$, by \eqref{rho}, and \eqref{u}, we have
\begin{equation}\begin{split}
D_1\leq\mid D_1 \mid=&\Big|\frac{\gamma}{\gamma - 1}\int_0^T \int_\Omega \psi(t)\Big(\overline{\rho\nabla(\rho^{\gamma-1})} - \rho \nabla(\rho^{\gamma-1})\Big)\cdot\overline{\mathbf{u}}\,\mathrm{d}x \,\mathrm{d}t\Big|\\
=&\Big|\gamma\int_0^T \int_\Omega\psi(t)\big[(\overline{\rho^{\gamma-1}\nabla\rho}-\rho^{\gamma-1}\nabla\rho)\cdot\overline{\mathbf{u}}
\big]\,\mathrm{d}x \,\mathrm{d}t\Big|\\
\leq&C\parallel\psi(t)\parallel_{L^{\infty}(0,T)}\int_0^T\parallel\overline{\mathbf{u}}\parallel_{L^q(\Omega)}
\parallel\overline{\rho^{\gamma-1}\nabla\rho}-\rho^{\gamma-1}\nabla\rho\parallel_{L^\frac{q}{q-1}(\Omega)}\mathrm{d}t\\
\leq&C\int_0^T\parallel\overline{\mathbf{u}}\parallel_{L^q(\Omega)}
\parallel\overline{\rho^{\gamma-1}\nabla\rho}-\rho^{\gamma-1}\nabla\rho\parallel_{L^2(\Omega)}\mathrm{d}t\\
\leq&C\parallel\overline{\mathbf{u}}\parallel_{L^p(0,T;L^q(\Omega))}
\parallel\overline{\rho^{\gamma-1}\nabla\rho}-\rho^{\gamma-1}\nabla\rho\parallel_{L^\frac{p}{p-1}(0,T;L^2(\Omega))}\rightarrow0,~~~as~~\epsilon\rightarrow0.
\end{split}\end{equation}

Using $\rho_t\in L^2(0,T;L^{\frac{2q}{q+2}}(\Omega))$ and $\rho\leq\tilde{\rho}$, we have $D_3$ goes to zero as $\epsilon$ tends to zero. Thus $R_{\epsilon}\to 0$, as $\epsilon \to 0$.

\medskip

Secondly, we prove $I_{\epsilon} +J_\epsilon \to 0$, as $\epsilon \to 0$. By the Gagliardo-Nirenberg inequality and \eqref{regu}, for any $\alpha_2\in(0,1),$ we obtain
\begin{equation}\label{regu of H}
  \textbf{H}\in L^{r_2}(0,T;L^{r_1}(\Omega)),
\end{equation}
where $\frac{1}{r_1}=\frac{1}{6}+\frac{\alpha_2}{3}, \frac{1}{r_2}=\frac{1-\alpha_2}{2}$.
In fact, for any $0<\alpha_2<1,$ we have
\begin{align*}
\parallel\textbf{H}\parallel^{r_2}_{L^{r_2}(0,T;L^{r_1}(\Omega))}
=& \int_0^T\parallel \textbf{H}\parallel_{L^{r_1}(\Omega)}^{r_2}\mathrm{d}t \\
\leq&\int_0^T\big(\parallel\textbf{H}\parallel^{\alpha_2}_{L^2(\Omega)}\parallel\textbf{H}\parallel^{1-\alpha_2}_{H^1(\Omega)}\big)^{r_2}\mathrm{d}t\\
\leq&C+\parallel \textbf{H}\parallel_{L^\infty(0,T;L^2(\Omega))}^{\alpha_2r_2}\int_0^T\parallel \nabla\textbf{H}\parallel_{L^2(\Omega)}^{r_2(1-\alpha_2)}\mathrm{d}t\\
\leq& C\int_0^T\parallel \nabla\textbf{H}\parallel_{L^2(\Omega)}^2\mathrm{d}t+C\\
\leq& C,
\end{align*}
where $\frac{\alpha_2}{2}+(\frac{1}{2}-\frac{1}{3})(1-\alpha_2)=\frac{1}{r_1},$ and $r_2(1-\alpha_2)=2.$

\medskip

Firstly, we prove $I_{\epsilon} =  I_1+I_2+I_3+I_4\rightarrow0,$ as $\epsilon \to 0$.

By virtue of the assumption that $\mathbf{u}\in L^p(0,T;L^q(\Omega))$ where $p,q$ will be determined later, and  H\"{o}lder inequality, Lemma \ref{limit} and \eqref{regu of H}, we get
\begin{align*}
I_{\epsilon}=&  I_1+I_2+I_3+I_4\leq \mid I_1\mid+\mid I_2\mid+\mid I_3\mid+\mid I_4\mid\\
=&\Big|\int_0^T \int_\Omega\psi(t)\Big(\mathrm{div}(\mathbf{H}\otimes\overline{\mathbf{H}})-\overline{\mathrm{div}(\mathbf{H}\otimes \mathbf{H})}\Big)\cdot\overline{\mathbf{u}}\,\mathrm{d}x \,\mathrm{d}t\Big|
+\Big|\frac{1}{2}\int_0^T \int_\Omega \psi(t)\overline{\mathbf{u}}\cdot\Big(\nabla\overline{(\mathbf{H}\cdot\mathbf{H})}
 - \nabla(\mathbf{H}\cdot \overline{\mathbf{H}})\Big)\,\mathrm{d}x \,\mathrm{d}t\Big|\\
&+\Big|\int_0^T \int_\Omega\psi(t)\overline{\mathbf{H}}\cdot\Big(\mathrm{div}(\mathbf{H}\otimes\overline{\mathbf{u}})
  -\overline{\mathrm{div}(\mathbf{H}\otimes\mathbf{u})}\Big)\,\mathrm{d}x \,\mathrm{d}t\Big|+\Big|\int_0^T \int_\Omega\psi(t)\overline{\mathbf{H}}\cdot\Big(\overline{\mathrm{div}(\mathbf{u}\otimes\mathbf{H})}
  -\mathrm{div}(\overline{\mathbf{u}}\otimes\mathbf{H})\Big)\,\mathrm{d}x \,\mathrm{d}t\Big|\\
  \leq& C\parallel\psi(t)\parallel_{L^\infty(0,T)}
  \int_0^T \parallel \textbf{u}\parallel_{L^q(\Omega)}\parallel\nabla \textbf{H}\parallel_{L^2(\Omega)}
  \parallel \textbf{H}\parallel_{L^\frac{2q}{q-2}(\Omega)}\mathrm{d}t\\
  \leq& C\parallel \textbf{u}\parallel_{L^p(0,T;L^q(\Omega))}\parallel\nabla \textbf{H}\parallel_{L^2(0,T;L^2(\Omega))}
  \parallel \textbf{H}\parallel_{L^\frac{2p}{p-2}(0,T;L^{r_1}(\Omega))}\\
  \leq& C\parallel \textbf{u}\parallel_{L^p(0,T;L^q(\Omega))}\parallel\nabla \textbf{H}\parallel_{L^2(0,T;L^2(\Omega))}
  \parallel \textbf{H}\parallel_{L^{r_2}(0,T;L^{r_1}(\Omega))},
\end{align*}
where $\frac{2q}{q-2}\leq r_1,$ and $\frac{2p}{p-2}\leq r_2.$
%$\frac{1}{q}\leq\frac{1}{3}-\frac{\alpha_1}{3},\frac{1}{p}\leq\frac{1}{2}.$
%$$\frac{3}{q}+\frac{2}{p}\leq1, q\geq6.$

Thanks to Lemma \ref{limit}, as $\epsilon$ tend to zero, we have
\begin{equation}\label{F}
  I_1+I_2+I_3+I_4\rightarrow0.
\end{equation}

Next we prove $J_\epsilon=J_1+J_2\to 0$, as $\epsilon \to 0$.\\
For $J_1$, we obtain
\begin{align*}
  \mid J_1\mid= & \mid\frac{1}{2}\int_0^T \int_\Omega\psi(t)(\overline{\mathbf{H}} - \mathbf{H})\cdot\Big((\overline{\mathbf{u}}\cdot\nabla)\mathbf{H}\Big)\mathrm{d}x\mathrm{d}t\mid\\
  \leq& C\parallel\psi(t)\parallel_{L^\infty(0,T)}
  \int_0^T \parallel \textbf{u}\parallel_{L^q(\Omega)}\parallel\nabla \textbf{H}\parallel_{L^2(\Omega)}
  \parallel \overline{\textbf{H}}-\textbf{H}\parallel_{L^\frac{2q}{q-2}(\Omega)}\mathrm{d}t\\
  \leq& C\parallel \textbf{u}\parallel_{L^p(0,T;L^q(\Omega))}\parallel\nabla \textbf{H}\parallel_{L^2(0,T;L^2(\Omega))}
  \parallel \overline{\textbf{H}}-\textbf{H}\parallel_{L^\frac{2p}{p-2}(0,T;L^{r_1}(\Omega))}\\
  \leq& C\parallel \textbf{u}\parallel_{L^p(0,T;L^q(\Omega))}\parallel\nabla \textbf{H}\parallel_{L^2(0,T;L^2(\Omega))}
  \parallel \overline{\textbf{H}}-\textbf{H}\parallel_{L^{r_2}(0,T;L^{r_1}(\Omega))}\rightarrow0,~~~as~~\epsilon\rightarrow0,
\end{align*}
where $\frac{2q}{q-2}\leq r_1,$ and $\frac{2p}{p-2}\leq r_2.$\\
For $J_2$, we have
\begin{align*}
\mid J_2\mid= & \mid\frac12\int_0^T \int_\Omega\psi(t)\mathbf{H}\cdot\Big((\overline{\mathbf{u}}\cdot\nabla)(\mathbf{H}-\overline{\mathbf{H}})\Big)\ dxdt\mid\\
  \leq& C\parallel\psi(t)\parallel_{L^\infty(0,T)}
  \int_0^T \parallel \textbf{u}\parallel_{L^q(\Omega)}\parallel \textbf{H}\parallel_{L^\frac{2q}{q-2}(\Omega)}
  \parallel \nabla(\overline{\textbf{H}}-\textbf{H})\parallel_{L^2(\Omega)}\ dt\\
  \leq& C\parallel \textbf{u}\parallel_{L^p(0,T;L^q(\Omega))}\parallel\textbf{H}\parallel_{L^\frac{2p}{p-2}(0,T;L^{r_1}(\Omega))}
  \parallel \nabla(\overline{\textbf{H}}-\textbf{H})\parallel_{L^2(0,T;L^2(\Omega))}\\
  \leq& C\parallel \textbf{u}\parallel_{L^p(0,T;L^q(\Omega))}\parallel\textbf{H}\parallel_{L^{r_2}(0,T;L^{r_1}(\Omega))}
  \parallel \nabla(\overline{\textbf{H}}-\textbf{H})\parallel_{L^2(0,T;L^2(\Omega))}\rightarrow0,~~~as~~\epsilon\rightarrow0,
\end{align*}
where $\frac{2q}{q-2}\leq r_1,$ and $\frac{2p}{p-2}\leq r_2.$

Thus, we have
\begin{equation}\label{J}
  J_\epsilon \to 0, ~as ~\epsilon \to 0.
\end{equation}
In fact, because $\frac{1}{r_1}=\frac{1}{6}+\frac{\alpha_1}{3}, \frac{1}{r_2}=\frac{1-\alpha_1}{2}$ is equivalent to
$\frac{3}{r_1}+\frac{2}{r_2}=\frac{3}{2}, r_1<6,$ and $\frac{1}{r_1}\leq\frac{1}{2}-\frac{1}{q},\frac{1}{r_2}\leq\frac{1}{2}-\frac{1}{p}$, we obtain
\begin{equation*}
  \frac{3}{2}=\frac{3}{r_1}+\frac{2}{r_2}\leq\frac{5}{2}-\bigg(\frac{3}{q}+\frac{2}{p}\bigg).
\end{equation*}
%and
%\begin{equation*}
%  \frac{1}{q} \leq\frac{1}{2}-\frac{1}{r_1}<\frac{1}{2}-\frac{1}{6}
%  <\frac{1}{3}.
%\end{equation*}
Thus we need $\mathbf{u}\in L^p(0,T;L^q(\Omega))$, where $\frac{3}{q}+\frac{2}{p}\leq1$ for any $p\ge4,\, q\geq6.$

We are ready to pass to the limits in \eqref{r40}. Let $\epsilon$ go to zero, and using \eqref{F}-\eqref{J}, what we have proved is that in the limit
\begin{equation}\begin{split}\label{test function}
&-\int_0^T \int_\Omega \psi_t\bigg(\frac{1}{2} \rho\mid\mathbf{u}\mid^2 +\frac{\rho^\gamma }{\gamma - 1} + \frac{1}{2}\mid\mathbf{H}\mid^2\bigg)\,\mathrm{d}x\,\mathrm{d}t\\
&+\int_0^T \int_\Omega\psi(t)\Big(\nu\mid\nabla\times\mathbf{H}\mid^2+ \mu\mid\nabla\mathbf{u}\mid^2+ (\lambda + \mu) \mid\mathrm{div}\mathbf{u}\mid^2\Big)\,\mathrm{d}x \,\mathrm{d}t=0
\end{split}\end{equation}
for any test function $ \psi\in\mathcal{D}(0, +\infty)$.\\

\textbf{Step 3.} Extending the result \eqref{test function} for $ \psi\in\mathcal{D}(-1, +\infty)$.

The final goal is to extend our result \eqref{test function} for the test function $ \psi\in\mathcal{D}(-1, +\infty)$. To this end, it is necessary for us to have the continuity of $\rho(t),(\sqrt{\rho}\mathbf{u})(t)$ and $\mathbf{H}$ in the strong topology at $t = 0$. Adopting a similar argument to that of \cite{Yu2017}, what we expected can be done.

Using $\sqrt{\rho}\mathbf{u}\in L^\infty(0,T;L^2(\Omega)),$ and (\ref{rho}), we have
\begin{equation*}
  \rho_t\in L^2(0,T;H^{-1}(\Omega)),\ and\ \nabla\rho\in L^2(0,T;L^2(\Omega)).
\end{equation*}
Hence
\begin{equation}\label{ST1}
  \rho\in C([0,T];L^2(\Omega)).
\end{equation}

By energy inequality \eqref{energy inequ}, we have $\sqrt{\rho}\mathbf {u}\in L^\infty(0,T;L^2(\Omega)),
\mathrm{div}\mathbf {u}\in L^2(0,T;L^2(\Omega))$, and $\rho\leq\tilde{\rho}.$ Recall $\eqref{2dbu-E1.1}_1$, we obtain
\begin{equation*}
  (\sqrt{\rho})_t=-\mathrm{div}(\sqrt{\rho}\mathbf {u})+\frac{1}{2}\sqrt{\rho}\mathrm{div} \mathbf{u}.
\end{equation*}
Hence, we deduce
\begin{equation*}
 (\sqrt{\rho})_t\in L^2(0,T; H^{-1}(\Omega)).
\end{equation*}
Meanwhile, due to $\nabla\sqrt{\rho}\in L^\infty(0,T;L^2(\Omega)),$ we get $\sqrt{\rho}\in C([0,T];L^2(\Omega))$.
More generally, in view of $\rho\in L^\infty(0,T;L^{\infty}(\Omega))$, we deduce
\begin{equation}\label{0}
  \rho\in C([0,T];L^s(\Omega)),
\end{equation}
where $1\leq s<\infty$.\\
On the other hand,
\begin{align*}
  &\frac{1}{2}ess\limsup_{t\rightarrow0}\int_\Omega\mid \sqrt{\rho}\mathbf{u}-\sqrt{\rho_0}\mathbf{u}_0\mid^2dx
  +ess\limsup_{t\rightarrow0}\int_\Omega\mid \mathbf{H}-\mathbf{H}_0\mid^2dx \\
  =& \underbrace{ess\limsup_{t\rightarrow0}\bigg[\int_\Omega\bigg(\frac{1}{2}\rho \mathbf{u}^2+\frac{\rho^\gamma}{\gamma-1}+\mid \mathbf{H}\mid^2\bigg)dx
  -\int_\Omega\bigg(\frac{1}{2}\rho_0 \mathbf{u}_0^2+\frac{\rho_0^\gamma}{\gamma-1}+\mid \mathbf{H}_0\mid^2\bigg)dx\bigg]}_{II_1}\\
  &+\underbrace{ess\limsup_{t\rightarrow0}\int_\Omega\mathbf{u}_0(\rho_0\mathbf{u}_0-\sqrt{\rho_0\rho}\mathbf{u})dx
  +2ess\limsup_{t\rightarrow0}\int_\Omega\mathbf{H}_0(\mathbf{H}_0-\mathbf{H})dx}_{II_2}\\
  &+\underbrace{ess\limsup_{t\rightarrow0}\int_\Omega\bigg(\frac{\rho_0^\gamma}{\gamma-1}-\frac{\rho^\gamma}{\gamma-1}\bigg)dx}_{II_3}.
\end{align*}
Using energy inequality \eqref{energy inequ}, we get $II_1\leq0.$

By \eqref{0}, we obtain
  $II_3=0.$\\
Hence, we have
\begin{align*}
  &\frac{1}{2}ess\limsup_{t\rightarrow0}\int_\Omega\mid \sqrt{\rho}\mathbf{u}-\sqrt{\rho_0}\mathbf{u}_0\mid^2dx
  +ess\limsup_{t\rightarrow0}\int_\Omega\mid \mathbf{H}-\mathbf{H}_0\mid^2dx \\
  %\leq&ess\limsup_{t\rightarrow0}\int_\Omega\sqrt{\rho_0}\mathbf{u}_0(\sqrt{\rho_0}\mathbf{u}_0-\sqrt{\rho}\mathbf{u})dx
  \leq&ess\limsup_{t\rightarrow0}\int_\Omega\mathbf{u}_0(\rho_0\mathbf{u}_0-\sqrt{\rho_0\rho}\mathbf{u})dx
  +2ess\limsup_{t\rightarrow0}\int_\Omega\mathbf{H}_0(\mathbf{H}_0-\mathbf{H})dx\\
  =&II_2^1+II_2^2.
\end{align*}
For $II_2^1,$ taking $\mathbf{u}_0\in L^\kappa,~\kappa>2$, it follows
\begin{align*}
  II_2^1= & ess\limsup_{t\rightarrow0}\int_\Omega\mathbf{u}_0(\rho_0\mathbf{u}_0-\rho\mathbf{u})
  +ess\limsup_{t\rightarrow0}\int_\Omega\mathbf{u}_0\sqrt{\rho}\mathbf{u}(\sqrt{\rho}-\sqrt{\rho_0})\\
  =&0,
\end{align*}
%\begin{align*}
%  I_2^1= & 2ess\limsup_{t\rightarrow0}\int_\Omega\frac{\sqrt{\rho_0}\mathbf{u}_0}{\sqrt{\rho}}(\sqrt{\rho\rho_0}\mathbf{u}_0-\rho_0\mathbf{u}_0)
%  +2ess\limsup_{t\rightarrow0}\int_\Omega\frac{\sqrt{\rho_0}\mathbf{u}_0}{\sqrt{\rho}}(\rho_0\mathbf{u}_0-\rho_0\mathbf{u}_0)\\
%  =&0,
%\end{align*}
where we used the fact $\rho\mathbf{u}\in C([0,T];L_{weak}^2(\Omega)), \sqrt{\rho}\mathbf{u}\in L^\infty(0,T;L^2(\Omega))$, and \eqref{0}.\\
\\
For $II_2^2,$ using $\mathbf{H}\in C([0,T];L_{weak}^2(\Omega))$, we get $II_2^2=0.$

So we get
\begin{equation*}
  \frac{1}{2}ess\limsup_{t\rightarrow0}\int_\Omega\mid \sqrt{\rho}\mathbf{u}-\sqrt{\rho_0}\mathbf{u}_0\mid^2dx
  +ess\limsup_{t\rightarrow0}\int_\Omega\mid \mathbf{H}-\mathbf{H}_0\mid^2dx \leq0,
\end{equation*}
which gives us
\begin{equation}\label{ST2}
  \sqrt{\rho}\mathbf{u}\in C([0,T];L^2(\Omega)), \ and\ \mathbf{H}\in C([0,T];L^2(\Omega)).
\end{equation}

%In the view of \eqref{ST1},\eqref{ST2} and Lebesgue  point theorem, we get
%\begin{equation}\label{lebsgue}
%  \lim_{\tau\rightarrow0}\frac{1}{\tau}\int_0^\tau\int_\Omega\bigg(\frac{1}{2}\rho \mathbf{u}^2+\frac{\rho^\gamma}{\gamma-1}+\frac{1}{2}\mid \mathbf{H}\mid^2\bigg)=\int_\Omega\bigg(\frac{1}{2}\rho_0 \mathbf{u}_0^2+\frac{\rho_0^\gamma}{\gamma-1}+\frac{1}{2}\mid \mathbf{H}_0\mid^2\bigg)dx
%\end{equation}

As in \cite{Wang2018}, for any given $t_0>0,$ we choose the text function $\psi_\tau(t)$ as below.
\begin{equation}\label{3.44}
\psi_{\tau}(t)=\left\{\begin{array}{cc}{0,} & {0 \leq t \leq \tau} \\

{\frac{t-\tau}{\epsilon},} & {\tau \leq t \leq \tau+\epsilon} \\

{1,} & {\tau+\epsilon \leq t \leq t_{0}} \\

{\frac{t_{0}+\epsilon-t}{\epsilon},} & {t_{0} \leq t \leq t_{0}+\epsilon} \\

{0,} & {t_{0}+\epsilon \leq t,}\end{array}\right.
\end{equation}
where positive $\tau$ and $\epsilon$ satisfying $\tau+\epsilon<t_0.$

Substituting (\ref{3.44}) into \eqref{test function}, we have
\begin{align}\label{new text function1}
  &\frac{1}{\epsilon}\int_{t_0}^{t_0+\epsilon}\int_\Omega\bigg(\frac{1}{2} \rho|\mathbf{u}|^2 +\frac{\rho^\gamma}{\gamma - 1}
  + \mid\mathbf{H}\mid^2\bigg)\,\mathrm{d}x \,\mathrm{d}t %\mid\mathbf{H}\mid^2\bigg)\,\mathrm{d}x\,\mathrm{d}t++\int_\tau^{\tau+\epsilon}\int_\Omega\frac{t-\tau}{\epsilon}\Big(\mu|\nabla\mathbf{u}|^2+ (\lambda + \mu) |\mathrm{div}\mathbf{u}|^2+\nu|\nabla\times\mathbf{H}|^2\Big)
 +\int_{\tau}^{t_0+\epsilon}\int_\Omega
  \psi_{\tau}(t)\Big(\mu|\nabla\mathbf{u}|^2+ (\lambda + \mu) |\mathrm{div}\mathbf{u}|^2+\nu|\nabla\times\mathbf{H}|^2\Big)\,\mathrm{d}x \,\mathrm{d}t\notag\\
 =& \frac{1}{\epsilon}\int_\tau^{\tau+\epsilon} \int_\Omega\bigg(\frac{1}{2} \rho|\mathbf{u}|^2
 +\frac{\rho^\gamma }{\gamma - 1}+ \mid\mathbf{H}\mid^2\bigg) \,\mathrm{d}x \,\mathrm{d}t .
  %+\int_{t_0}^{t_0+\epsilon}\int_\Omega\frac{t_0-t}{\epsilon}
 %\Big(\mu|\nabla\mathbf{u}|^2+ (\lambda + \mu) |\mathrm{div}\mathbf{u}|^2+\nu|\nabla\times\mathbf{H}|^2\Big)\,\mathrm{d}x \,\mathrm{d}t.
\end{align}
We deal with the second item on the left hand side
\begin{align*}
  &\int_{\tau}^{t_0+\epsilon}\int_\Omega
  \psi_{\tau}(t)\Big(\mu|\nabla\mathbf{u}|^2+ (\lambda + \mu) |\mathrm{div}\mathbf{u}|^2+\nu|\nabla\times\mathbf{H}|^2\Big)\,\mathrm{d}x \,\mathrm{d}t\\
  = & \int_{\tau}^{t_0+\epsilon}\int_\Omega
  \Big(\mu|\nabla\mathbf{u}|^2+ (\lambda + \mu) |\mathrm{div}\mathbf{u}|^2+\nu|\nabla\times\mathbf{H}|^2\Big)\,\mathrm{d}x \,\mathrm{d}t\\
  &+\int_{\tau}^{\tau+\epsilon}\int_\Omega
  \Big(\frac{t-\tau}{\epsilon}-1\Big)\Big(\mu|\nabla\mathbf{u}|^2+ (\lambda + \mu) |\mathrm{div}\mathbf{u}|^2+\nu|\nabla\times\mathbf{H}|^2\Big)\,\mathrm{d}x \,\mathrm{d}t\\
  &+\frac{1}{\epsilon}\int_{t_0}^{t_0+\epsilon}\int_\Omega
  (t_0-t)\Big(\mu|\nabla\mathbf{u}|^2+ (\lambda + \mu) |\mathrm{div}\mathbf{u}|^2+\nu|\nabla\times\mathbf{H}|^2\Big)\,\mathrm{d}x \,\mathrm{d}t.
\end{align*}
Passing to the limit as $\epsilon \to 0$ in \eqref{new text function1}, and using \eqref{ST1},\eqref{ST2} and Lebesgue theorem, one deduces
\begin{align}\label{new text function2}
  &\int_\Omega\bigg(\frac{1}{2} \rho|\mathbf{u}|^2
 +\frac{\rho^\gamma }{\gamma - 1}+ \mid\mathbf{H}\mid^2\bigg)(t_0) \,\mathrm{d}x+\int_\tau^{t_0}\int_\Omega
  \Big(\mu|\nabla\mathbf{u}|^2+ (\lambda + \mu) |\mathrm{div}\mathbf{u}|^2+\nu|\nabla\times\mathbf{H}|^2\Big)\,\mathrm{d}x \,\mathrm{d}t\notag\\
 =& \int_\Omega\bigg(\frac{1}{2} \rho|\mathbf{u}|^2 +\frac{\rho^\gamma}{\gamma - 1}
  + \mid\mathbf{H}\mid^2\bigg)(\tau)\,\mathrm{d}x,
\end{align}
where
\begin{align*}
  &\frac{1}{\epsilon}\int_{\tau}^{\tau+\epsilon}\int_\Omega
  (t-\tau)\Big(\mu|\nabla\mathbf{u}|^2+ (\lambda + \mu) |\mathrm{div}\mathbf{u}|^2+\nu|\nabla\times\mathbf{H}|^2\Big)\,\mathrm{d}x \,\mathrm{d}t\\
  \leq& \int_{\tau}^{\tau+\epsilon}\int_\Omega
  \Big(\mu|\nabla\mathbf{u}|^2+ (\lambda + \mu) |\mathrm{div}\mathbf{u}|^2+\nu|\nabla\times\mathbf{H}|^2\Big)\,\mathrm{d}x \,\mathrm{d}t \to 0,
\end{align*}
and similarly
\begin{align*}
  &\frac{1}{\epsilon}\int_{t_0}^{t_0+\epsilon}\int_\Omega
  (t_0-t)\Big(\mu|\nabla\mathbf{u}|^2+ (\lambda + \mu) |\mathrm{div}\mathbf{u}|^2+\nu|\nabla\times\mathbf{H}|^2\Big)\,\mathrm{d}x \,\mathrm{d}t\\
  \geq& \int_{t_0}^{t_0+\epsilon}\int_\Omega
  \Big(\mu|\nabla\mathbf{u}|^2+ (\lambda + \mu) |\mathrm{div}\mathbf{u}|^2+\nu|\nabla\times\mathbf{H}|^2\Big)\,\mathrm{d}x \,\mathrm{d}t \to 0,
\end{align*}
as $\epsilon \to 0.$
Finally setting $\tau\rightarrow0$ in \eqref{new text function2}, from \eqref{ST2} we get \eqref{energy inequ-1}.

\section*{Acknowledgements} The authors would like to thank Professors Changjiang Zhu and Huanyao Wen for some helpful suggestion.
Y. Chen is supported by the National Natural Science
Foundation of China (No. 11601160), by Science and
Technology Program of Guangzhou (No. 201707010221). M. Zhang is supported by the National Natural Science
Foundation of China (No. 11701185).


\addcontentsline{toc}{section}{\\References}

\begin{thebibliography}{00}
%% \bibitem{label}
%% Text of bibliographic item
%\bibitem{Alfven1}H. Alfv\'{e}n. Existence of electromagnetic-hydrodynamic waves, Nature, 150 (1942), 405-406.
%\bibitem{Alfven}H. Alfv\'{e}n. Cosmical Electrodynamics. Oxford University Press, Oxford, 1950.
\bibitem{Bardos2018}
C.Bardos, E.S.  Titi. Onsager's conjecture for the incompressible Euler equations in bounded domains. Arch Ration Mech Anal.,
228 (2018), 197-207.
\bibitem{Buckmaster2015}T. Buckmaster, C. Lellis, P. Isett, Jr.L. Szekelyhidi.  Anomalous dissipation for 1/5-H\"{o}lder Euler flows. Ann. Math., 182 (2015), 127-172.
\bibitem{Buckmaster2016}T. Buckmaster, C. Lellis, L. Sz\'{e}kelyhidi. Dissipative Euler flows with Onsager-critical spatial regularity. Comm. Pure Appl. Math., 69 (2016), 1613-1670.
\bibitem{Cabannes}H. Cabannes. Theoretical Magnetofluiddynamics. Academic Press, New York, 1970.
\bibitem{Wang2018} R.M. Chen, Z.L. Liang, D.H. Wang and R.Z. Xu. Energy equality in compressible fluids with physical boundaries,	arXiv:1808.06089v2, 2018.
\bibitem{CW1}G.Q. Chen, D.H. Wang. Global solutions of nonlinear magnetohydrodynamics with large initial data. J. Differential Equations, 182 (2002), 344-376.
\bibitem{CW2}G.Q. Chen, D.H. Wang. Existence and continuous dependence of large solutions for the magnetohydrodynamic equations. Z. Angew. Math. Phys., 54 (2003), 608-632.
\bibitem{Cheskidov}A. Cheskidov, P. Constantin, S. Friedlander and R. Shvydkoy. Energy conservation and Onsager's conjecture for the Euler equations. Nonlinearity, 21 (2008), 1233-1252.
\bibitem{Constantin}P. Constantin, E. Weinan and E.S. Titi. Onsager's conjecture on the energy conservation for solutions of Euler's equation. Commun. Math. Phys., 165 (1994), 207-209.
%\bibitem {Du2014}L.L. Du, Y.F. Wang. Blow-up criterion for 3-dimensional compressibel Navier-Stokes equations involving velocity divergence. Commun. Math. Sci., 12 (2014), 1427-1435.
%\bibitem{DuvautandLions} G. Duvaut, J.L. Lions. Inequation en theremoelasticite et magnetohydrodynamique. Arch. Rational Mech. Anal., 46 (1972), 241-279.
\bibitem{Eyink}Eyink, L. Gregory. Energy dissipation without viscosity in ideal hydrodynamics I. Fourier analysis and local energy transfer. International Symposium on Physical Design, 1994.
\bibitem{Fan}J. Fan, W. Yu. Strong solution to the compressible magnetohydrodynamic equations with vacuum. Nonlinear Analysis: Real World Application, 10 (2009), 392-409.
\bibitem{Feireisl}E. Feireisl. Dynamics of viscous compressible fluids. Oxford University Press, Oxford, 2004.
\bibitem{Feir2}E. Feireisl, A. Novotn\'{y}, H. Petzeltov\'{a}. On the existence of globally defined weak solutions to the Navier-Stokes equations. J. Math. Fluid Mech., 3 (2001), 358-392.
\bibitem{hu2008}X.P. Hu, D.H. Wang. Compactness of weak solutions to the three-dimensional compressible magnetohydrodynamic equations. J. Differential Equations, 245 (2008), 2176-2198.
\bibitem{Hu2008}X.P. Hu, D.H. Wang. Global solutions to the three-dimensional full compressible magnetohydrodynamic flows. Comm. Math. Phys., 283 (2008), 255-284.
\bibitem{Wang2010} X.P. Hu, D.H. Wang. Global existence and large-time behavior of solutions to the three-dimensional equations of compressible magnetohydrodynamic flows. Arch. Ration. Mech. Anal., 197 (2010), 203-238.
%\bibitem{Isett}P. Isett. A proof of Onsager's conjecture. Preprint, 2016.
\bibitem{Kawashima}S. Kawashima, M. Okada. Smooth global solutions for the one-dimensional equations
in magnetohydrodynamics. Proc. Japan Acad. Ser. A Math. Sci., 58 (1982), 384-387.
\bibitem{Kulikovskiy}A.G. Kulikovskiy, G.A. Lyubimov. Magnetohydrodynamics. Addison-Wesley,  Massachusetts, 1965.
\bibitem{Laudau}L.D. Laudau, E.M. Lifshitz. The electrodynamics of Continuous Media, 2nd ed., Pergamon, New York, 1984.
\bibitem{Li2013} H.L. Li, X.Y. Xu and J.W. Zhang. Global classical solutions to 3D compressible magnetohydrodynamic equations with large oscillations and vacuum. SIAM J. Math. Anal., 45 (2013), 1356-1387.
\bibitem{Lions}P.L. Lions. Mathematical topics in fluid mechanics. Vol. 2. Compressible models. Oxford University Press, New York, 1998.
\bibitem{Onsager}L. Onsager. Statistical hydrodynamics. Nuovo Cimento, 6 (1949), 279-287.
\bibitem{Serrin}J. Serrin. The initial value problem for the Navier-Stokes equations. Nonlinear Problems. Proceedings of the Symposium, Madison, Wisconsin, 1962. University of Wisconsin Press, Madison, Wisconsin, 69-98, 1963.
\bibitem{Shinbrot}M. Shinbrot. The energy equation for the Navier-Stokes system. SIAM J. Math. Anal., 5 (1974), 948-954.
\bibitem{Vasseur}A.F. Vasseur, C. Yu. Existence of global weak solutions for 3D degenerate compressible Navier-Stokes equations. Invent. Math., 206 (2016), 935-974.
\bibitem{Volpert}A.I. Vol\'{p}ert, S.I. Hudjaev. The Cauchy problem for composite systems of nonlinear differential equations. Mat. Sb., 87 (1972), 504-528.
\bibitem{Yu2017} C. Yu. Energy conservation for the weak solutions of the compressible Navier-Stokes equations. Arch. Ration. Mech. Anal., 225 (2017), 1073-1087.
\end{thebibliography}
\end{document}